\documentclass[12pt]{amsart}
\usepackage{amssymb}
\DeclareMathAlphabet{\mathantt}{OT1}{antt}{li}{it}
\DeclareMathAlphabet{\mathpzc}{OT1}{pzc}{m}{it}

\usepackage{tikz}
\usepackage[colorlinks=true,linkcolor=blue,urlcolor=blue,citecolor=blue, pagebackref]{hyperref}
\usepackage[alphabetic]{amsrefs}
\usepackage{soul}
\usepackage[margin=1.0in]{geometry}
\usepackage{xypic, verbatim,amscd,color}

\newcommand\wtu{\widetilde{U}}

\newcommand\starr{\bullet}

\newcommand\dlog{\operatorname{dlog}}

\newcommand\op{\operatorname}

\DeclareMathOperator{\srhom}{\mathcal{R}\text{\kern -1pt {hom}}\,}

\newcommand\frg{\mathfrak{g}}
\newcommand\frh{\mathfrak{h}}

\newcommand\wtd{\widetilde{\Delta}}
\newcommand\wtp{\widetilde{P}}

\newcommand\mI{\mathcal{I}}

\newcommand\mf{\mathcal{F}}

\newcommand\mk{\mathcal{K}}

\newcommand\wD{\widetilde{\Delta}}

\newcommand{\Poincare}{Poincar\'e}

\newcommand{\Ccech}{\v{C}ech }

\newcommand\tensor{\otimes}
\newcommand\ml{\mathcal{L}}

\newcommand\wP{\widetilde{P}}

\newcommand\mg{\mathcal{G}}

\newcommand\im{\operatorname{im}}

\newcommand{\leto}[1]{\stackrel{#1}{\to}}

\newtheorem{theorem}[equation]{Theorem}
\newtheorem{remark}[equation]{Remark}

\newtheorem{corollary}[equation]{Corollary}
\newtheorem{question}[equation]{Question}
\newtheorem{proposition}[equation]{Proposition}

\newtheorem{lemma}[equation]{Lemma}

\newtheorem{definition}[equation]{Definition}
\newtheorem{defi}[equation]{Definition}
\newtheorem{observation}[equation]{Observation}

\begin{document}
\title{Hyperplane arrangements and tensor product invariants}
\author{P. ~Belkale, P.~Brosnan and S.~Mukhopadhyay}
\maketitle

\begin{abstract}
  In the first part of this paper, we consider, in the context of an
  arbitrary hyperplane arrangement, the map from compactly
  supported cohomology to the usual cohomology of a local system. A
  formula (i.e., an explicit algebraic de Rham representative) for a
  generalized version of this map is obtained.

These results are applied in the second part to invariant theory: Schechtman and Varchenko connect invariant theoretic objects to the cohomology of local systems on complements of hyperplane arrangements. The first part of this paper is then used, following and completing arguments of Looijenga, to determine the image of invariants in cohomology. In suitable cases (e.g., corresponding to positive integral levels), the space of invariants is shown to acquire a mixed Hodge structure over a cyclotomic field. We investigate the Hodge filtration on the  space of invariants, and characterize the subspace of conformal blocks in Hodge theoretic terms.

\end{abstract}
\section{Introduction}
Let $W$ be a $M$-dimensional  complex vector space, and consider an arbitrary weighted hyperplane arrangement $(W,\mathcal{C},a)$ in $W$. This consists of the following data:
\begin{enumerate}
\item[(a)] A collection $\mathcal{C}$ of hyperplanes in $W$. We will assume that we are given polynomials $f_1,\dots,f_r$  on $W$ of degree one such that the hyperplanes in the collection $\mathcal{C}$ are the zero loci $Z(f_1),\dots,Z(f_r)$. These are subsets of $W$, which need not pass through the origin.
\item[(b)] A  vector $a=(a_1,\dots,a_r)\in\mathbb{C}^r$, which we will call a weighting of $\mathcal{C}$. The weight of  the hyperplane $Z(f_i)$ is $a_i$.
\end{enumerate}
Let $U=W-Z(f_1.f_2\cdots f_r)$.
The differential form corresponding to the weighting is
$$
\eta=\sum_{i=1}^r a_i \dlog f_i =\sum_{i=1}^r a_i\frac{df_i}{f_i} \in H^0(U,\Omega^1).
$$
Now set $\nabla=d+\eta:\mathcal{O}_{U^{\mathrm{an}}}\to \Omega_{U^{\mathrm{an}}}$ and write $\mathcal{L}(a)$ for the kernel of $\nabla$ which is a rank one local system on $U$.

\begin{defi}
The Aomoto subalgebra of the weighted hyperplane arrangement $(W,\mathcal{C},a)$ is the subalgebra $A^{\starr}(U)$ of $\Omega^{\starr}(U)$ generated over $\mathbb{C}$ by the forms $\dlog f_i$, for $i=1,\dots, r$. It is equipped with a differential $d\omega:=\eta\wedge \omega$. We call the pair $(A^{\starr}(U), \eta \wedge)$, the Aomoto complex.
\end{defi}
\subsection{The case of small weights}
To motivate the discussion, let us assume first that the weights $a_i$ are sufficiently small, or more generally, assume that the absolute values $|\sum_{i=1}^{r} \epsilon_i a_i|\in \Bbb{C}-\Bbb{Z}_{>0}$ for all choices of $\epsilon_i\in \{0,1\}$
 (we need the weights and their negatives to be  satisfy the conditions (Mon) in \cite{ESV}). In this case it is known by results of Esnault-Schechtman-Viehweg \cite{ESV} that the $i$-th cohomology group of the Aomoto complex represents $H^i(U,\mathcal{L}(a))$. This generalizes results of E. Brieskorn \cite{Brieskorn}.
One consequence of this is that $H^i(U,\mathcal{L}(a))=0$ for $i>M$.

Now consider the natural map
\begin{equation}\label{naturale}
H^i_c(U,\mathcal{L}(a))\to H^i(U,\mathcal{L}(a)).
\end{equation}
The compactly supported cohomology $H^i_c(U,\mathcal{L}(a))$ is  \Poincare-Verdier dual to
$H^{2M-i}(U,\mathcal{L}(a)^*)$ and $\mathcal{L}(a)^*=\mathcal{L}(-a)$. By \cite{ESV} applied to $\mathcal{L}(-a)$, we get that $H^{2M-i}(U,\mathcal{L}(-a))=0$ for $i<M$. Therefore $H^i_c(U,\mathcal{L}(a))$ also vanishes for $i<M$.

The map \eqref{naturale} is therefore interesting only when $i=M$. In this case it is induced by an element 
\begin{equation}\label{doppioo}
\Sigma\in H^{M}(U,\mathcal{L}(a)^*)\otimes H^M(U,\mathcal{L}(a))= H^{2M}(U\times U,\mathcal{L}(a)\boxtimes\mathcal{L}(-a)).
\end{equation}

Now $H^{2M}(U\times U,\mathcal{L}(a)\boxtimes\mathcal{L}(-a))$
is computed by the Aomoto complex of a hyperplane arrangement in $W\times W$: The hyperplanes are of the form $H\times W$ and $W\times H$ with weights $a_i$ and $-a_i$ for $H=Z(f_i)$. It therefore makes sense to ask for a formula for an explicit element in $A^{2M}(U\times U)$ which gives rise to the element $\Sigma$ in \eqref{doppioo}.

\begin{theorem}\label{maino}
Let
\begin{equation}\label{diagonalSIntro}
S=\sum_{1\leq i_1<\dots <i_{M}\leq r}\prod_{s=1}^{M} a_{i_s}\dlog f^{(1)}_{i_s}\dlog f^{(2)}_{i_s}
\end{equation}
(here $f^{(1)}_i(x,y)=f_i(x)$ and $f^{(2)}_i(x,y)=f_i(y)$). There exists a non zero $c\in\Bbb{C}$ such that element $cS\in A^{2M}(U\times U)$ represents the cohomology class  $\Sigma\in H^{2M}(U\times U,\mathcal{L}(a)\boxtimes\mathcal{L}(-a))$.
\end{theorem}
\begin{remark}
It is easy to see using \cite[Page 651]{GH} (see Lemma \ref{GHresult} below)  that $c=(2\pi \sqrt{-1})^{-M}$. The exact value of $c$ is not needed in this paper.
\end{remark}
The Aomoto complex computes topological cohomological groups even when the weights $a$ are not small, as shown by Looijenga \cite[Proposition 4.2]{L1}. We refer the reader to Section \ref{Lgen} for more details. The groups $H^i(U,\mathcal{L}(a))$ need to be replaced by hypercohomologies of  complexes on suitable compactifications of $U$, and we have a generalization of the map \eqref{naturale} again in this set-up. Our main result, Theorem \ref{main}, gives an explicit form which represents (up to a constant multiplicative scalar) this generalization of the mapping \eqref{naturale} (see also the map \eqref{zhGC}). In Remark \ref{saturdayy} we note that Theorem \ref{maino} holds for arbitrary weights.

\subsection{The log form $S$}
The map $A^M(U)^*\to A^M(U)$ induced by the element $S$ (in Theorem \ref{maino}) appears in the work of Schechtman and Varchenko \cite[Theorem 2.4 and  Lemma 3.2.5]{SV}. There is an entire collection of such  forms, our $S$ corresponds to $``{S^M}"$ in loc. cit.; all of these show up in the proof of Theorem \ref{maino} in roughly the following way:
\begin{itemize}

\item The element $\Sigma$ is a suitable cohomology class of the diagonal in $U\times U$, and should therefore vanish on the open subsets $U\times U- \{(x,y)\in U\times U\mid F(x)=F(y)\}$, where $F$ is an arbitrary linear form on $W$.

\item Therefore one should have (for Theorem \ref{maino} to hold)  $S$ exact in the Aomoto complex that computes the cohomology of the local system on the complements of the hyperplane arrangement given by deleting hyperplanes of the $\{(x,y)\in W\times W\mid F(x)=F(y)\}$ with weight $0$ from $U\times U$. Recall that $U\times U$ is the complement of the product weighted arrangement with weights $a_i$ on the hyperplane $Z(f_{i}^{(1)})$ and $-a_i$ on $Z(f_{i}^{(2)})$ for $1\leq i\leq r$.  We refer the reader to Corollary \ref{grundlegendSpeciale} for more details.

\item The proof that $S$ is exact features the other  forms $S^{(b)}$ (one needs to also consider intersections of such open subsets). The above argument is carried out in \Ccech cohomology, which leads to a full proof of Theorem \ref{maino}.

\end{itemize}
\subsection{The case of arbitrary weights}\label{Lgen}
We give a description of Looijenga's results from  \cite[Proposition 4.2]{L1} sufficient for the statement of a generalization of \eqref{naturale}.

Let $P$ be any smooth projective compactification of $U$, with
$P-U=\cup_{\alpha} E_{\alpha}$ a divisor with normal crossings.  Let
$V=P-\cup'_{\alpha}E_{\alpha}$, where the union is restricted to
$\alpha$ such that $a_{\alpha}=\op{Res}_{E_{\alpha}}\eta$ is not a
strictly positive integer. Similarly let
$V'=P-\cup'_{\alpha}E_{\alpha}$, where the union is restricted to
$\alpha$ such that $a_{\alpha}=\op{Res}_{E_{\alpha}}\eta$, the residue
of $\eta$ on $E_{\alpha}$, is not an integer which is $\geq 0$. Note
that $V'\supseteq V$.

Let $q:U\to V'$ and $j:U\to V$ denote the inclusion map. The
cohomology of the Aomoto complex $H^{\starr}(A^{\starr}(U),\eta)$ is
equal to the cohomology $H^{\starr}(V,j_{!}\mathcal{L}(a))$. The map
\eqref{naturale} is replaced by the natural map (which is non-zero
only for $\starr=M$)
\begin{equation}\label{bedEq}
H^{\starr}(V',q_{!}\mathcal{L}(a))\to H^{\starr}(V,j_{!}\mathcal{L}(a)).
\end{equation}
The groups $H^M(V',q_{!}\mathcal{L}(a))$ are dual to Aomoto cohomology groups (for the weight vector $-a$), and this set-up generalizes the case of small weights (if the weights are small then $V=U$, and $H^{M}(V',q_{!}\mathcal{L}(a))=H^{M}_c(U,\ml(a))$, $H^{M}(V,j_{!}\mathcal{L}(a))=H^{M}(U,\ml(a))$). Theorem \ref{main}, which is stated in this context computes \eqref{bedEq} and generalizes Theorem \ref{maino} for arbirtrary weights. In fact Theorem \ref{maino} is also true for all choices of weights, see Remark \ref{saturdayy}.
\subsection{Applications to invariant theory}\label{IntroIn}

Consider a finite dimensional simple Lie algebra $\frg$ with a fixed Cartan decomposition and let $R$ denote the set of positive simple roots. Let $(\ ,\ )$ be a normalized Killing form on $\frh$, the Lie algebra of the Cartan subgroup, such that $(\theta,\theta)=2$ where $\theta\in \frh^*$ is the highest root (identifying $\frh$ and $\frh^*$ using the Killing form).
For a dominant integral weight $\lambda\in\frh^*$, let $V_{\lambda}$ denote the corresponding irreducible representation.

Now suppose that we are given an $n$-tuple $\vec{\lambda}=(\lambda_1,\dots,\lambda_n)$ of dominant integral weights of $\frg$. The space of coinvariants
$$\Bbb{A}(\vec{\lambda})=(V_{\lambda_1}\tensor V_{\lambda_2}\tensor\dots\tensor V_{\lambda_n})_\frg=\frac{(V_{\lambda_1}\tensor V_{\lambda_2}\tensor\dots\tensor V_{\lambda_n})}{\frg (V_{\lambda_1}\tensor V_{\lambda_2}\tensor\dots\tensor V_{\lambda_n})}$$
is a fundamental object of invariant theory. Note that the space of invariants maps isomorphically onto the space of coinvariants under the obvious map. It is easy to see that $\Bbb{A}(\vec{\lambda})$ is zero if $\sum\lambda_i$ is not a positive sum of simple roots, and we will assume that this is indeed the case and write
$$\mu=\sum_{p=1}^r n_p \alpha_p,\ \ n_p\in\Bbb{Z}_{\geq 0}, \ \mbox{and} \ \alpha_p \in R.$$
  Fix a  map $\beta:[M]=\{1,\dots,M\} \to R$, so that
  \begin{equation}\label{hurricane}
  \mu=\sum_{b=1}^M \beta(b),\ \mbox{and} \  M=\sum_{p=1}^r n_p.
  \end{equation}

  Fix a point $\vec{z}=(z_1,\dots,z_n)$ in the configuration space of
  $n$ distinct points on $\Bbb{A}^1$. Let $W=\Bbb{C}^M$. The
  coordinate variables of $W$ will be denoted by $t_1,\dots,t_M$. We
  will consider the variable $t_b$ to be colored by the simple root
  $\beta(b)$. Let $\kappa$ be an arbitrary non-zero complex number and
  consider the weighted hyperplane arrangement $(W,\mathcal{C},a)$ in
  $W$ given by the following collection of hyperplanes, and their
  attached weights:
\begin{enumerate}
\item For $i\in [1,n]$ and $b\in [1,M]$, the hyperplane $t_b-z_i=0$, with weight $\displaystyle \frac{(\lambda_i,\beta(b))}{\kappa}$.
\item $b,c\in [1,M]$ with $b<c$,  the hyperplane $t_b-t_c=0$, with weight $\displaystyle-\frac{(\beta(b),\beta(c))}{\kappa}$.
\end{enumerate}
Let $U$ be the complement of the above hyperplane arrangement in $W$ as before. The corresponding $\eta\in A^1(U)$ is
$$\eta=\frac{1}{\kappa}\biggl(\sum_{i=1}^n\sum_{b=1}^M (\lambda_j,\beta(b))\frac{d(t_b-z_i)}{t_b-z_i} - \sum_{b,c\in [1,M], b<c} (\beta(b),\beta(c))\frac{d(t_b-t_c)}{t_b-t_c}\biggr) .$$

The basic connection between invariant theory and the topology of hyperplane arrangements arises from an injective map  (see Section \ref{SVM} for a description ) constructed by Schechtman and Varchenko \cite{SV}:
\begin{equation}\label{SVmap}
\Omega^{SV}:{V}(\vec{\lambda})^*_0=(V_{\lambda_1}\tensor V_{\lambda_2}\tensor\dots\tensor V_{\lambda_n})^*_0\to A^M(U).
\end{equation}

Here ${V}(\vec{\lambda})^*_0$ is the zero weight space (for $\frh$)  in the dual of ${V}(\vec{\lambda})= (V_{\lambda_1}\tensor V_{\lambda_2}\tensor\dots\tensor V_{\lambda_n})$. 
Note that
$\Bbb{A}(\vec{\lambda})^*$ is a subspace of ${V}(\vec{\lambda})^*_0$.
\begin{remark}
The map in \cite[Theorem 6.16.2]{SV} is for Verma modules for the corresponding Lie algebra without Serre relations, and is in this context an isomorphism. This isomorphism is also described explictly in \cite[Proposition 5.7]{BM}. The representations $V_{\lambda}$ considered here are quotients of these Verma modules, and hence we get an injective map in \eqref{SVmap}.

\end{remark}
Now since $U$ is of dimension $M$ , we get that $A^M(U)$ surjects onto $H^M(A^{\starr}(U), \eta\wedge)$. Using unitarity results on conformal blocks \cite{Ram,Bel}, we show injectivity in the following:
\begin{proposition}\label{injective}
The induced mapping $\Bbb{A}(\vec{\lambda})^{*}\to H^M(A^{\starr}(U),\eta \wedge)$ is injective.
\end{proposition}

\begin{remark}
Proposition \ref{injective} above appears in \cite[Lemma 3.3]{L2}. The sketched proof of this proposition in \cite{L2}is inadequate: The proof in \cite{L2} is via \cite[Proposition 2.18]{L2}, which is an injection into a term in the Aomoto complex. But injection in a quotient of the top degree term in this complex is more subtle (indeed, the proof in our paper uses extension to compactifications).
\end{remark}

 There is a natural symmetric group acting on $H^M(A^{\starr}(U),\eta \wedge)$:
\begin{equation}\label{symmetrie}
\Sigma_M=\{\sigma\in S_M\mid \beta(\sigma(b))=\beta(b), b=1,\dots,M\}.
\end{equation}
We consider $\Sigma_M$ as the  color preserving symmetric group acting on variables $t_1,\dots,t_M$. The map  $\Bbb{A}(\vec{\lambda})^*\to H^M(A^{\starr}(U),\eta \wedge)$ has its image in $H^M(A^{\starr}(U),\eta \wedge)^{\chi}$ where $\chi$ is the sign character on $\Sigma_M$ induced by the symmetric group $S_M$. Therefore one obtains an injective map
\begin{equation}\label{beforeclass}
\Bbb{A}(\vec{\lambda})^*\hookrightarrow H^M(A^{\starr}(U),\eta\wedge)^{\chi}.
\end{equation}
Note that as in Section \ref{Lgen},  we may write $H^M(A^{\starr}(U),\eta\wedge)^{\chi}=H^M(V,j_{!}\mathcal{L}(a))^{\chi}$. In the case of small weights, this becomes $H^M(A^{\starr}(U),\eta \wedge)^{\chi}=H^M(U,\mathcal{L}(a))^{\chi}$.
So the right hand side of \eqref{beforeclass} has a topological interpretation. Although this will not play a role in this paper, we note that the induced mapping
\begin{equation}\label{LD2}
\Bbb{A}(\vec{\lambda})^*\to H^M(V,j_{!}\mathcal{L}(a))^{\chi}
\end{equation}
is flat for connections as $\vec{z}=(z_1,\dots,z_n)$ varies in the configuration space of $n$-distinct points on $\Bbb{A}^1$. Here the left hand side of \eqref{LD2} has the Knizhnik-Zamolodchikov (KZ) connection, and the right hand side the Gauss-Manin connection \cite[Lemma 6.6.3]{SV}, \cite[Lemma 3.9, Lemma 3.10]{L2}.

\subsection{The image of the injective map (14)}
Looijenga's strategy for  determining this image is the following: For simplicity of exposition we will assume, in the introduction, that  $\kappa$ is sufficiently large in absolute value so that the weights above are small (this assumption is eventually dropped). Therefore,
$H^M(A^{\starr}(U),\eta \wedge)\leto{\sim}H^M(U,\mathcal{L}(a))$  (also similarly for $-a$)
and hence under the assumption of smallness of weights, there is an injection
\begin{equation}\label{LD4}
\Bbb{A}(\vec{\lambda})^*\hookrightarrow H^M(U,\mathcal{L}(a))^{\chi}.
\end{equation}
 Consider the map \eqref{beforeclass} for the dual weights $\lambda_1^*,\dots,\lambda_n^*$, and $-\kappa$ for the value of $\kappa$. Since $\lambda^*=-w_0\lambda$, where $w_0$ is the longest element in the Weyl group, we can write
    $$\sum_{i=1}^n \lambda_i =\sum_{b=1}^M\gamma(b).$$
    where $\gamma(b)=-w_0\beta(b)$ are simple positive roots. Therefore we may use the same vector space and variables $t_1,\dots,t_M$, just the colors change for the dual. The weights are now
\begin{enumerate}
\item For $i\in [1,n]$ and $b\in [1,M]$, the hyperplane $t_b-z_i=0$, with weight $\displaystyle \frac{(-w_0\lambda_i,-w_0\beta(b))}{-\kappa}$ which is also equal to $\displaystyle\frac{-(\lambda_i, \beta(b))}{\kappa}$.
\item $b,c\in [1,M]$ with $b<c$,  the hyperplane $t_b-t_c=0$, with weight $\displaystyle-\frac{(-w_0\beta(b),-w_0\beta(c))}{-\kappa}$ which is also equal to $\displaystyle \frac{-(\beta(b),\beta(c))}{\kappa}$.
\end{enumerate}Equalities of weights  above hold because the Cartan-Killing form is invariant under the Weyl group. These weights are negatives of the weights assigned for $\lambda_1,\dots,\lambda_n$ and $\kappa$.  We therefore have an injection
$$\Bbb{A}(\vec{\lambda}^*)^*\hookrightarrow H^M(A^{\starr}(U),-\eta\wedge)^{\chi}=H^M(U,\mathcal{L}(-a))^{\chi}.$$
Dualizing, we find a surjection ($\chi$ is a sign character, hence self dual)
\begin{equation}\label{LD1}
H^M_c(U,\mathcal{L}(a))^{\chi}\twoheadrightarrow  \Bbb{A}(\vec{\lambda}^*).
\end{equation}
Now, there is a canonical isomorphism from invariants to coinvariants
\begin{equation}\label{canaan}
\Bbb{A}(\vec{\lambda})^*\leto{\sim}\Bbb{A}(\vec{\lambda}^*).
 \end{equation}
 Putting the maps \eqref{LD1} and \eqref{LD2} together with the inverse of \eqref{canaan}, we get
\begin{equation}\label{LD3}
H^M_c(U,\mathcal{L}(a))^{\chi}\twoheadrightarrow  \Bbb{A}(\vec{\lambda}^*)\leto{\sim}\Bbb{A}(\vec{\lambda})^*\hookrightarrow H^M(U,\mathcal{L}(a))^{\chi}.
\end{equation}

Therefore (following Looijenga \cite[Theorem 3.7]{L2}), we see that the image of the map \eqref{LD2} equals the image of the composite $H^M_c(U,\mathcal{L}(a))^{\chi}\to  H^M(U,\mathcal{L}(a))^{\chi}$ in \eqref{LD3}.

Now Looijenga assumes, without proof, the following compatibility property:
\begin{itemize}
 \item the composite map in \eqref{LD3} is the natural
$H^M_c(U,\mathcal{L}(a))^{\chi}\to  H^M(U,\mathcal{L}(a))^{\chi}$, induced by topology (see \cite[Proof of Theorem 3.7, page 33]{L2}).
\end{itemize}
 and  concludes the following (actually a generalized form,
valid for arbitrary $\kappa\neq 0$ is proved in \cite{L2}, this result is Theorem \ref{egregium} in our paper):
\begin{theorem}\label{secondmain}
Assume $|\kappa|$ is sufficiently large. The image of invariants in the topological cohomological groups, i.e., the image of the injective map \eqref{LD4}, coincides with the image of the map induced by topology
$H^M_c(U,\mathcal{L}(a))^{\chi}\to  H^M(U,\mathcal{L}(a))^{\chi}$. Therefore,
$$\Bbb{A}(\vec{\lambda})^*=\operatorname{Image}\ H^M_c(U,\mathcal{L}(a))^{\chi}\to  H^M(U,\mathcal{L}(a))^{\chi}.$$
\end{theorem}
\begin{remark}
Looijenga's assumption is subtle, and needs a justification since the map \eqref{LD3} factors through a representation
theoretic duality \eqref{canaan}, and also uses \Poincare-Verdier duality and the maps \eqref{SVmap}. Therefore Looijenga's proof assumes a compatibility property between topological and representation theoretic dualities, as well as a compatibility with Schechtman-Varchenko maps \eqref{SVmap}.
\end{remark}

The results of the first part of the paper were motivated by the problem of proving this interesting compatibility property. In Theorems \ref{main} and \ref{maino}, which are valid for arbitrary weighed hyperplane arrangements, a formula for the (topological) $\Sigma_M$ equivariant  mapping $H^M_c(U,\mathcal{L}(a))\to  H^M(U,\mathcal{L}(a))$
(and a generalization) is obtained. By this we mean a de Rham representative for the corresponding element (which will be $\Sigma_M$-invariant for the action on $U\times U$) of the space \eqref{doppioo}.

This formula is compared with a formula for the actual composite \eqref{LD3} which is obtained using the work of Schechtman and Varchenko \cite[Theorem 6.6]{SV}. These formulas coincide (up to a non-zero scalar), and one proves the assumption implicit in  Looijenga's proof of \cite[Theorem 3.7]{L2}.

\subsection{Rational weights and the Hodge theory of the KZ system}
In Section \ref{mathew}, we show that our results imply that $\Bbb{A}(\vec{\lambda})^*$ carries a mixed Hodge structure over a cyclotomic field extension of $\Bbb{Q}$ if $\kappa$ is an integer (or even a rational number). If $\kappa=\ell +g^*$ where $\ell$ is a positive integer and $g^*$ the dual Coxeter number of $\frg$, then using results in \cite{BM}, we show that the $F^M$ (Hodge) part intersected with $W_M$ (weight) part of  $\Bbb{A}(\vec{\lambda})^*$ coincides with the space of conformal blocks at level $\ell$, for $\frg$ classical or the Lie algebra of the exceptional group $G_2$. While we have some expectations for the Hodge filtration (see Question \ref{plausible}), the weight filtration on $\Bbb{A}(\vec{\lambda})^*$ remains mysterious, and one could ask for a formula for the Hodge numbers $h^{p,q}$.

In Section \ref{PS} we give an example where the mixed Hodge structure on $\Bbb{A}(\vec{\lambda})^*$ is not pure, by showing that the monodromy of the KZ system is not semisimple.

\section{The Aomoto complex}
The aim of this section is to make explicit what the Aomoto complex represents.

\subsection{Resolution of singularities}\label{TonyS}
Let $(X,D)$ be a pair of a smooth variety $X$ and a divisor $D$. There exists \cite{BVP} a canonical resolution of singularities of $(X,D)$: a birational projective morphism  $\pi:\widetilde{X}\to X$ with the following properties
\begin{enumerate}
\item[(a)] $\widetilde{X}$ is smooth.
\item[(b)]  $\widetilde{D}=\pi^{-1}(D)$ is a divisor with simple normal crossings.
\item[(c)]  $\pi:\widetilde{X}-\widetilde{D}\to X-D$ is an isomorphism.
\item[(d)] For any $p\in D$ such that $D$ is simple normal crossings at $p$, the map $\pi$ is an isomorphism over a neighborhood of $p$.
\item[(e)] Automorphisms of the pair $(X,D)$, extend to the pair $(\widetilde{X},\widetilde{D})$.
\end{enumerate}
We will have occasion to use the above result when $D$ is not locally a hyperplane arrangement. Property (d) is used in an essential manner in this paper (Section \ref{classo}); Property (e) can be avoided, but leads to a more satisfying picture.

\subsection{Constuctible complexes}
If $X$ is a complex algebraic variety, we let
$D^b_c(X)$ denote the bounded derived category of sheaves of $\Bbb{C}$-vector spaces on $X$
with constructible cohomology. We refer the reader to  \cite[Section 4.5, Page 111]{Hota} for more details on the definition and notation.

\subsection{Compactifications}\label{flashback1}
Let $P$ be any smooth projective compactification of $U$, with $P-U=\cup_{\alpha} E_{\alpha}$ a divisor with normal crossings.
The higher cohomology $H^j(P,\Omega_{P}^i(\log E))$, $j>0$ vanishes by \cite[Section 2]{ESV}, and $H^0(P,\Omega_{P}^i(\log E))$ is the space of log forms $A^i(U)$. The global hypercohomology of the complex $(\Omega^{\starr}_{{P}}(\log E),d+\eta)$ coincides with the cohomology of the Aomoto complex. In Lemma
\ref{beanT} and \ref{stringbeans}, we recall known results on the stalks of the Aomoto complex. These are used to give a topological characterization of the Aomoto complex in Lemma \ref{beans}. We regard $\Omega^{i}_P(\log E)$ as sheaves in the analytic topology unless explicitly mentioned.

The following Lemma is a part of \cite[Proposition 3.13]{DiffEq}:
\begin{lemma}\label{stringbeans}
Suppose $p\in P$. Assume that  $a_{\beta}$ is not a strictly positive integer for any $E_{\beta}$ passing through $p$. Then the local hypercohomology at $p$ of the complex $(\Omega^{\starr}_{{P}}(\log E),d+\eta)$ is isomorphic to the hypercohomology at $p$ of $Rk_{*}\ml(a)$ where $k:U\to P$ (via the canonical map to $Rk_{*}\ml(a)$ obtained by adjunction).
\end{lemma}

The next Lemma (as well as its proof) is a generalization of Lemma~\ref{stringbeans}.

\begin{lemma}\label{beanT}
Suppose $p\in E_{\beta}\subseteq P$. Assume $a_{\beta}$ is not an integer which is $\leq 0$. Then the stalk of the hypercohomology at $p$ of the complex $(\Omega^{\starr}_{{P}}(\log E),d+\eta)$ is zero.
\end{lemma}
Note that in the above statement we allow $p$ to be in the intersection of several $E_{\alpha}$; For some (but not all) of these, $a_{\alpha}$ could be an integer $\leq 0$.

\begin{proof}
  We can replace $P$ by the open polydisc $\Delta^n$, and assume that
  $E$ is contained in the union of coordinate hyperplanes with
  $E_{\beta}$ the hyperplane given by $z_1=0$.  In fact, we write the
  proof for the case where $E$ is exactly the union of the
  hyperplanes, leaving to the reader the obvious modifications
  necessary for the general case.  We can write
  $\eta=\sum a_i \dlog z_i$, and, by our hypothesis,
  $a_1$ is not a negative integer.

Our aim is to construct a retraction $\delta:\Omega^k_{\Delta^n}(\log E)))_0
\to \Omega^{k-1}_{\Delta^n}(\log E))_0$ of the stalk of the complex
$(\Omega^{\bullet}_{\Delta^n}(E),d+\eta)$ at $0$.
To this end, let $\theta$ denote the endomorphism of $\Omega^k_{\Delta^n}(\log E)))_0$ obtained by contraction with the vector field $\displaystyle
z_1\frac{\partial}{\partial z_1}$.
Then, for a multi-index $I=(i_1,\ldots, i_k)$ with $1\leq i_1 < i_2<\cdots < i_k\leq n$,
$$
  \theta (\dlog z_{i_1}\wedge \cdots \wedge \dlog z_{i_k})=
\begin{cases}
  \dlog z_{i_2}\wedge \cdots \wedge \dlog z_{i_k}, & i_1=1\\
                                                0, & \text{else.}
\end{cases}
$$
For a vector $m=(m_1,\ldots, m_n)\in\mathbb{Z}_{\geq 0}^n$, set
$z^m:=z_1^{m_1}z_2^{m_2}\ldots z_n^{m_n}$.  Then a general form
$\alpha\in\Omega^k_{\Delta^n}(\log E)_0$ can be written as a convergent
power series $$\alpha=\sum_{m,I}c_{m,I} z^m \dlog z_{i_1}\cdots \dlog z_{i_k}$$
where here the $c_{m,I}$ are in $\mathbb{C}$ and we drop the exterior product
symbol for brevity
of notation.
We set
$$
\delta\alpha:=\sum_{m,I} \frac{c_{m,I}}{m_1+a_1} z^m \theta(\dlog z_{i_1}\cdots \dlog z_{i_k}).
$$
Note that the convergence of the power series representation for $\alpha$
(in an open polydisk containing $0$) immediately implies the convergence of the above power series representation for $\delta\alpha$ (in the same polydisk).

Now, we claim that
\begin{equation}
  (d+\eta)\delta\alpha +\delta (d+\eta)\alpha=\alpha.
  \label{htpy}
\end{equation}
The operators  $d$, $d+\eta$ and $\delta$ are all $\mathbb{C}$-linear
and continuous with respect to the
$\mathfrak{m}$-adic topology where $\mathfrak{m}$ is the maximal ideal
in $\mathcal{O}_{\Delta^n,0}$. So it suffices to check~\eqref{htpy}
on monomials of the form  $\alpha=z^m\dlog z_{i_1}\dots \dlog z_{i_k}$.
For this, the computation is as follows:
\begin{itemize}
\item[(a)] If $i_1\neq 1$, set $\beta=z_2^{m_2}\dots z_n^{m_n} \dlog
  z_{i_1}\dots \dlog z_{i_k}$. Clearly
  $(d+\eta)\delta (\alpha)=0$. Now,
  $\delta (d+\eta)(\alpha)=\delta (m_1z_1^{m_1} \dlog
  z_1\wedge \beta +a_1 z_1^{m_1} \dlog z_1\wedge \beta)=\alpha$.

\item[(b)] If $i_1=1$, set $\beta=z_2^{m_2}\dots z_n^{m_n} \widehat{\dlog
    z_{1}}\dots \dlog z_{i_k}$ then
  $$(d+\eta)\delta\alpha=
\frac{1}{m_1+\alpha_1}\left( z_1^{m_1} d\beta +\eta' \beta
  +(m_1+a_1)\alpha\right),$$
where here $\eta'=\sum_{i>1}a_i \dlog
  z_i$. Now the terms in $(d+\eta)(\alpha)$ which contain $\dlog
  z_1$ are of the form $-(\dlog z_1) z_1^{m_1} d\beta-(\dlog
  z_1)z_1^{m_1}\eta'\beta$. The desired equality follows.
\end{itemize}
\end{proof}

\begin{definition}
Write $T$ for the stratification of $\Delta^n$ determined by the
coordinate hyperplanes.  So,
the strata $T_{\alpha}$ are indexed by subsets $\alpha\subset\{1,\ldots, n\}$,
and $(z_1,\ldots, z_n)\in T_{\alpha}$
if $\alpha=\{i:z_i\neq 0\}$.
Write $D^b_T(\Delta^n)$ for the full subcategory of
  $D^b_c(\Delta^n)$ constisting of complexes of sheaves whose
  cohomology is constructible with respect to $T$.  In other words, a complex
  $\mathcal{F}\in D^b_c(\Delta^n)$ is in $D^b_T(\Delta^n)$ if its
  cohomology sheaves are locally constant on the strata.
Similarly, if $U$ is any locally closed union of strata of $T$,
we write $D^b_T(U)$ of the full subcategory of $D^b_c(U)$ consisting
of complexes whose cohomology is constructible with respect to the stratification on $U$ induced by $T$.
\end{definition}

\begin{observation}\label{obsone}
For $\mathcal{F}\in D^b_T(\Delta^n)$ and
$i\in \mathbb{Z}$,
$H^i(\Delta^n,\mathcal{F})= H^i(\mathcal{F}_0)$.
In other words the cohomology of $\mathcal{F}$ on $\Delta^n$ is the
same as the cohomology of the stalk at the origin.
\end{observation}

\begin{proof}
Write $\Delta_{\epsilon}$ for the disk of radius $\epsilon$ centered at the
origin in $\mathbb{C}$. Then
$$H^i(\mathcal{F}_0)=\varinjlim_{0<\epsilon<1} H^i(\Delta_{\epsilon}^n,\mathcal{F}).$$
But all the terms in the colimit are isomorphic because $\mathcal{F}$ is
constructible with respect to the stratification by the hyperplanes.
\end{proof}

\begin{corollary}\label{cortoobsone}
Write $i:\{0\}\to \Delta^n$ for the inclusion.
If $\mathcal{F}\in D^b_T(\Delta^n)$ and $i^*\mathcal{F}$ is quasi-isomorphic
to $0$, then $H^k(\Delta^n,\mathcal{F})=0$ for all integers $k$.
\end{corollary}

\begin{proof}
  By Observation~\ref{obsone}, the cohomology of the complex
$i^*\mathcal{F}$ is isomorphic to $H^*(\Delta^n,\mathcal{F})$.
\end{proof}

\begin{corollary}\label{ctoo}
Suppose $U$ is a locally closed subset of $\Delta^n$ which is a union
 strata of $T$ not containing the origin.  Write $j:U\to \Delta^n$
for the inclusion of $U$.  Then, for any, complex $\mathcal{F}\in D^b_T(U)$,
$H^i(\Delta^n,j_!\mathcal{F})=0$.
\end{corollary}

\begin{proof}
  The functor $j_!$ is exact and, for any sheaf $\mathcal{G}$ on $U$,
  the stalk of $j_!\mathcal{G}$ at the origin vanishes.  So the
  complex $i^*j_!\mathcal{F}$ is quasi-isomorphic to $0$ (with
  $i:\{0\} \to\Delta^n$ the inclusion).  The result follows from
  Corollary~\ref{cortoobsone}.
\end{proof}

\begin{lemma}\label{OnOrOff}
Consider $\eta= \sum_{i=1}^n a_i \log z_i$ on $\Delta^n\times \Delta^m$. Let, for $0\leq i\leq n$,
$$U_i=(\Delta^*)^{i}\times \Delta^{n-i}\times \Delta^m=\{(z_1,\dots,z_{n+m})\mid\  |z_a|<1,z_b\neq 0,\forall a,\forall 1\leq b\leq i\},\ U_0=\Delta^n\times \Delta^m $$
($\Delta^m$ corresponds to the set of variables that do not appear in the boundary divisor).
We have maps $j_i:U_i\to U_{i-1}$.
Consider an object $\mf=Rj_{1,?}Rj_{2,?}\dots Rj_{n,?}\ml(a)$, where the $?$ are possibly different choices from $\{!,*\}$. Suppose at least one of the $?$ is a $!$. Then, the local cohomology of $\mf$ at $0\in \Delta^n\times \Delta^m$ is zero.
\end{lemma}

\begin{proof}
  Assume that $i$ is the smallest index such that the $?$ in
  $Rj_{i,?}$ is a $!$. The statement is obvious if $i=1$, and
  therefore we assume $i>1$. Consider
  $k:U_{i-1}\to U_0=\Delta^n\times \Delta^m$, and
  $\mg= Rj_{i,!}Rj_{i+1,?}\dots Rj_{n,?}\ml(a)$ on $U_{i-1}$. We need
  to show that $H^{\starr}(U_0,Rk_*\mg)=H^{\starr}(U_{i-1},\mg)=0$.

  Let $p:U_{i-1}\to (\Delta^*)^{i-1}$ be the projection. It suffices
  to show that $Rp_* \mg=0$, which is local on $(\Delta^*)^{i-1}$. We
  may therefore assume $a_1=\dots=a_{i-1}=0$, and compute fiber wise
  along $p$. We need to show that
  $H^{\starr}(\Delta^{n-i+1}\times \Delta^m,j_!\mk)=0$, where
  $j:\Delta^*\times \Delta^{n-i}\times \Delta^m\to
  \Delta^{n-i+1}\times\Delta^m$
  and $\mk$ is a suitable complex of sheaves. This vanishing follows
  from Corollary~\ref{ctoo}.
\end{proof}

Let $V=P-\cup'_{\alpha}E_{\alpha}$, where $a_{\alpha}$ is not a strictly positive integer. Let $j:U\to V$, $k:V\to P$  denote the inclusion maps. The following is a result of Looijenga:
\begin{lemma}\label{beans}The complex $\Omega_{\eta}^{\starr}=(\Omega^{\starr}_{{P}}(\log E),d+\eta)$ equals $Rk_*j_{!}\mathcal{L}(a)$ as objects in $D^b_c(P)$.
\end{lemma}

\begin{proof} %
(\cite[Proposition 4.2]{L1})
There is a canonical isomorphism $\mathcal{L}(a)\to \Omega_{\eta}^{\starr}$ on $U$. By adjunction this gives rise to a map
$j_{!}\mathcal{L}(a)\to \Omega_{\eta}^{\starr}$ on $V$ which is verified to be an quasi-isomorphism using Lemma \ref{beanT}. We get isomorphisms
$$Rk_{*}j_{!}\mathcal{L}(a)\to Rk_*\Omega_{\eta}^{\starr}.$$

Now there is by adjunction a canonical map  $\Omega_{\eta}^{\starr}\to Rk_*\Omega_{\eta}^{\starr}$. We claim this is an isomorphism, and hence complete the proof of Lemma \ref{beans}. Once this claim is proved, the quasi-isomorphism from $Rk_{*}j_{!}\mathcal{L}(a)$ to $\Omega_{\eta}^{\starr}$ is obtained in two steps: the map $Rk_{*}j_{!}\mathcal{L}(a)\leto{\sim} Rk_*\Omega_{\eta}^{\starr},$ composed with the inverse of
$\Omega_{\eta}^{\starr}\leto{\sim} Rk_*\Omega_{\eta}^{\starr}.$

We prove this claim by comparing local cohomologies at points $p\in P$. Let $p\in P-V$. If $p$ does not lie on an $E_{\beta}$ with $a_{\beta}$ a positive integer, the desired isomorphism follows from Lemma \ref{stringbeans}.

Now assume that $p\in E_{\beta}$ and $a_{\beta}$ a positive
integer. Then the local cohomology at $p$ of $\Omega_{\eta}^{\starr}$
is zero by Lemma \ref{beanT}. Therefore to conclude the proof we only
have to verify that the stalk at $p$ of $Rk_{*}j_{!}\mathcal{L}(a)$ is
zero. This follows from Lemma~\ref{OnOrOff}.
\end{proof}

\subsection{A variation}\label{financial}
We could have taken $\hat{V}=P-\cup E_{\alpha}$ all $\alpha$ such that $a_{\alpha}$ is either a strictly positive integer, or a non-integer. Let $j':U\to \hat{V}$ and $k':\hat{V}\to P$.
Then $\Omega_{\eta}^{\starr}$ is represented by $k'_{!}Rj'_*\ml(a)$. This comes about by combining the isomorphisms
$$k'_{!}\Omega_{\eta}^{\starr}\leto{\sim} k'_{!}Rj'_*\ml(a) \ \ \mbox{and} \ \ k'_{!}\Omega_{\eta}^{\starr}\leto{\sim}\Omega_{\eta}^{\starr}.$$

\subsection{}\label{mostgeneralform}
Let $E_1,\dots,E_n$ be an enumeration of the irreducible components of $P-U$. The index set of $\alpha$ is therefore $\{1,\dots,n\}$. Let $U_0=P$, $U_1=P-E_1$, $U_2=P-E_1\cup E_2$, etc, and $U=U_n=P=E_1\cup E_2\cup\dots\cup E_n$. Let $j_i:U_i\to U_{i-1}$.

Color the divisors $E_{\alpha}$ by four colors (the divisor is colored, not the points on it!):
\begin{itemize}
\item If $a_{\alpha}$ is not an integer, then color the divisor {green},
\item If $a_{\alpha}$ is a positive integer then color the divisor {white},
\item If $a_{\alpha}$ is a negative integer, color the divisor black,
\item If $a_{\alpha}$ equals zero, color the divisor blue.
\end{itemize}

Now form an extension of the sheaf $\ml(a)$ on $U$ to all of $P$ as follows: The extension is
\begin{equation}\label{bigex}
Rj_{1,?}Rj_{2,?}\dots Rj_{n,?}\ml(a).
\end{equation}
Here the $?$ in $Rj_{i,?}$ is $!$ if the color on $E_i$ is white, and $*$ if black or blue, and either $!$ or $*$ if the color is green (both produce the same answer). Note that $j_{!}$ is exact and $j_!=Rj_{!}$.

We claim that the resulting object in the derived category is independent of the ordering of divisors using standard adjunction properties (see Remark \ref{adjunction} below). To prove this we consider $V=P-\cup' E_{\alpha}$, with $E_{\alpha}$  either green or white. All of the sheaves produced have zero local cohomology at points of $\cup' E_{\alpha}$, therefore all the sheaves are extension by zero from $V$. Restricted to $V$, all are derived lower star extensions  from $U$ which obviously commute.

\begin{defi}\label{underline}
Let $\underline{\ml}(a)=(\Omega^{\starr}_{{P}}(\log E),d+\eta)$ as an object in $D^b_c(P)$. This object is canonically quasi-isomorphic to any of the elements \eqref{bigex} as above.
\end{defi}

\begin{remark}\label{adjunction}If  $k: V\subset  P$ is any open subset and $F\in D^b_c(V)$ and $G\in D^b_c(P)$ then we have adjunctions \cite[Proposition 2.6.4 and 3.1.12]{KS}
$\operatorname{Hom}_P(k_{!}F,G)= \operatorname{Hom}_V(F,G_V)$ and
$\operatorname{Hom}_P(G,Rk_*F)= \operatorname{Hom}_V(G_V,F)$, where $G_V$ is the pull back of $G$ to the open subset $V$.
\end{remark}

\subsection{Independence from choices of compactifications}\label{officehours}
Suppose ${P}'$ is another compactification of $U$  with $P'-U=\cup_{\beta} {E}'_{\beta}$ a divisor with normal crossings. Assume that there is a map  $\pi:{P}'\to P$ which is identity over $U$.

We know that $(\Omega^{\starr}_{{P}}(\log E),d+\eta)$ represents $\underline{\ml}(a)$. There is an evident morphism
$$(\Omega^{\starr}_{{P}}(\log E),d+\eta)\to (\pi_*\Omega^{\starr}_{{P'}}(\log E'),d+\eta)\to R\pi_* (\Omega^{\starr}_{{P'}}(\log E'),d+\eta)$$
and hence a morphism
\begin{equation}\label{canonR}
 \underline{\ml}(a)\to R\pi_*\underline{\ml'}(a).
\end{equation}
This morphism obviously identifies  $H^i(P, \underline{\ml}(a))$ and $H^i(P', \underline{\ml'}(a))$ compatibly with  Aomoto cohomology.

By \cite[Lemma VI.3.4]{FC}, for any $i$, $\Omega^i_{{P'}}(\log E')$ is acyclic for the functor $\pi_*$ (i.e., higher direct images vanish), and $\pi_*\Omega^i_{{P'}}(\log E')= \Omega^i_{{P}}(\log E)$. Therefore,
\begin{proposition}\label{seventy}
 The morphism \eqref{canonR} is a quasi-isomorphism.
\end{proposition}

\subsection{A topological construction of the morphism (34)}
In the rest of the section, we construct an inverse to \eqref{canonR}, and construct \eqref{canonR} topologically.

Let
$$V'=P'-\cup' E'_{\beta},\ \  V= P-\cup' E_{\alpha}$$
with the union restricted to $\beta$ and $\alpha$ with $a'_{\beta}$ and $a_{\alpha}$ in the set $\Bbb{C}-\{1,2,3,\dots\}$. Let $j:U\to V$ and $j':U\to V'$.

We construct a natural morphism
\begin{equation}\label{canon}
R\pi_*\underline{\ml'}(a)\to \underline{\ml}(a).
\end{equation}
It suffices to construct such a morphism over $V$, by adjunction properties of $Rk_*$ where $k:V\to P$. By Lemma \ref{overman} below, it is lower shriek extension through out in $V$ and $\pi^{-1}(V)$, and the map \eqref{canon} is evident over $V$(compare stalks of both sides on $V-U$, and show they are both zero by proper base change).
\begin{lemma}\label{overman}
$\pi^{-1}(V)\subseteq V'$
\end{lemma}
\begin{proof}
We claim that for any divisor $E'_{\alpha'}$ which intersects $V'$, $a_{\alpha'}$ is a positive integer. Now $\pi(E_{\alpha'})$ has its generic point in $V$. Let $p'$ be a generic point of $E'_{\alpha}$ and $p=\pi(p')$. Assume we have coordinate systems
$u_1,\dots,u_M$ on $P'$ and $z_1,\dots,z_M$ near $p'$ and $p$ respectively so that $u_1=0$ is $E'_{\alpha_1}$, and $z_j,1\leq j\leq s$ are the divisors $E_{\bullet}$  passing though $p$. Therefore $z_j$ for $1\leq j\leq s$ pull back to functions divisible by $u_1$ (in a neighborhood of $p'$). Also the zeros of $z_j=0$ for $1\leq j\leq s$ are contained in $E'_{\alpha'}$. Hence we may write for $1\leq i\leq s$
$$z_i=u_1^{m_i}f(z')$$ with $f(u)$ invertible near $p'$ and $m_i>0$.  The residue of the pullback of such $\dlog z_i$ is $m_i$, and $\eta$ pulls back to a form with residue along $E'_{\alpha'}$ given by a linear combination $\sum_{i=1}^s m_i a_i$
where $a_i$ is the residue of $\eta$ along $z_i=0$ ($1\leq i\leq s$). This linear combination is clearly a positive integer as desired.
\end{proof}

There is another way of getting canonical morphisms in this picture: Let $$\hat{V}'=P'-\cup' E'_{\beta},\ \  \hat{V}= P-\cup' E_{\alpha}$$
with the union restricted to $\beta$ and $\alpha$ with $a'_{\beta}$ and $a_{\alpha}$ in the set $\Bbb{C}-\{0,-1,-2,-3,\dots\}$. Let $j:U\to \hat{V}$ and $j':U\to \hat{V}'$. Using this set-up, we may create a morphism which goes in a direction opposite to \eqref{canon}
\begin{equation}\label{canonRR}
 \underline{\ml}(a)\to R\pi_*\underline{\ml'}(a).
\end{equation}
Using adjunction properties  of $j_!$ (see Section \ref{financial}), it suffices to construct such a morphism over $\hat{V}$. Over $\hat{V}$, the residues over divisors in $P'$ are sums of elements in $\{0,-1,-2,-3,\dots\}$ (by using the same argument as in Lemma \ref{overman}), and hence we extend using $Rk_*$ from $U$ in $P'$, therefore we are again done using
adjunction.

It is also easy to see (by adjunction properties again) that \eqref{canonRR} coincides with \eqref{canonR}.

\begin{proposition}\label{trail}
The morphism \eqref{canon} is the inverse of \eqref{canonR}.
\end{proposition}
\begin{proof}
It suffices  to show that the composition
\begin{equation}\label{canonRRR}
 \underline{\ml}(a)\to R\pi_*\underline{\ml'}(a)\to\underline{\ml}(a)
\end{equation}
is an quasi-isomorphism. By various adjunction properties,
$$\operatorname{Hom}_P(\underline{\ml}(a),\underline{\ml}(a))=
\operatorname{Hom}_U({\ml}(a),{\ml}(a)).$$
therefore any homomorphism $\underline{\ml}(a)\to \underline{\ml}(a)$ is determined by its restriction to $U$ and the proposition follows.
\end{proof}

\section{The main morphism}
\begin{defi}
Let $X$ be a  complex algebraic variety.
For an element $F\in D^b_c(X)$, the Verdier dual $DF$ is defined as follows
$$DF= {{ {{\srhom}}}}(F, a_X^{!}\Bbb{C}),$$
 where $a_X:X\to \operatorname{Spec}(\Bbb{C})$.
\end{defi}
Recall that for any open inclusion  $j:V\to X$, and $F\in D^b_c(V)$, $D(j_{!}F)=Rj_*D(F)$ and $D(Rj_{*}F)=j_{!}D(F)$. Moreover if $\ml$ is a local system on $V$,  $D(\ml)[-2M]=\ml^*$ where $M$ is the dimension of $V$. Using the description \eqref{bigex} of $\underline{\ml}(-a)$, we see that the Verdier dual ${D}(\underline{\ml}(-a))[-2M]$ is almost the same as $\underline{\ml}(a)$, the only difference is extension over the blue divisors (divisors $E_{\alpha}$ such that $a_{\alpha}=0)$): In  $\underline{\ml}(a)$, the extension is $Rj_*$ and in ${D}(\underline{\ml}(-a))$, it is $j_!$. Therefore, in the notation of Section \ref{mostgeneralform} (with colors on divisors for $a$)
\begin{equation}\label{bigex2}
{D}(\underline{\ml}(-a))[-2M]=Rj_{1,?}Rj_{2,?}\dots Rj_{n,?}\ml(a).
\end{equation}
Here the $?$ in $Rj_{i,?}$ is $!$ if the color on $E_i$ is white or blue, and $*$ if black, and either $!$ or $*$ if the color is green (both produce the same answer). Therefore there is an obvious map
\begin{equation}\label{zh}
{D}(\underline{\ml}(-a))[- 2M]\to \underline{\ml}(a).
\end{equation}
Recall that $M$ is the complex dimension of $P$. The following lemma shows that up to scale there is only possible morphism between the objects that appear in \eqref{zh}.
\begin{lemma}
\begin{equation}\label{dyson2}
\operatorname{Hom}_P\bigl({D}(\underline{\ml}(-a))[-2M], \underline{\ml}(a)\bigr)=\Bbb{C}.
\end{equation}

\end{lemma}
\begin{proof}
We first note the following: $X$ be a  complex algebraic variety, $j: V\to X$ an open inclusion, and  $F,G\in D^b_c(V)$,
Then $j_!$ and $Rj_*$ are fully faithful:
$\operatorname{Hom}_X(j_!F,j_!G)=\operatorname{Hom}_V(F,G)$ and similarly for $Rj_*$. These properties follow from adjunction
(see Remark \ref{adjunction}). In addition we also have (again by adjunction)
$\operatorname{Hom}_X(j_!F,Rj_*G)=\operatorname{Hom}_V(F,G)$.

By the description \eqref{bigex2} and \eqref{bigex} of ${D}(\underline{\ml}(-a))$ and $\underline{\ml}(a)$ respectively, we
now get
\begin{equation}\label{dyson1}
\operatorname{Hom}_P({D}\bigl(\underline{\ml}(-a)\bigr)[-2M], \underline{\ml}(a))=\operatorname{Hom}_U\bigl({D}(\underline{\ml}(-a))[-2M],\underline{\ml}(a)\bigr).
\end{equation}

But ${D}(\underline{\ml}(-a))[-2M]$ restricted to $U$ equals the local system $\ml(a)$, and hence the last space in
\eqref{dyson1} equals $$\operatorname{Hom}_U(\ml(a),{\ml}(a)\bigr)=\Bbb{C}.$$
\end{proof}
\begin{remark}\label{todayy}
By \cite[Appendix A]{EV},
$${D}(\underline{\ml}(-a))[-2M]\leto{\sim}(\Omega^{\starr}_{{P}}(\log E)(-E),d+\eta)$$
and the map \eqref{zh} is the natural map,
\begin{equation}\label{zh1}
(\Omega^{\starr}_{{P}}(\log E)(-E),d+\eta)\to (\Omega^{\starr}_{{P}}(\log E),d+\eta).
\end{equation}
\end{remark}

\subsection{Small weights}
If the weights $a$ are sufficiently small, then $H^M(P,\underline{\ml}(a))=H^M(U,\ml(a))$ because there $j_{!}$ extensions are not used. Similarly, $H^M(P,D(\underline{\ml}(-a))[-2M])= H^M_c(U,\ml(a))$.

\subsection{The map (39) in global cohomology}\label{afternoon}
\begin{equation}\label{zhGC}
H^i(P,{D}(\underline{\ml}(-a))[-2M]) = H^{2M-i}(P, \underline{\ml}(-a))^*\to H^i(P,\underline{\ml}(a)).
\end{equation}
This map is non-zero only for $i=M$: the RHS vanishes for $i>M$ since it is computed by Aomoto cohomology, and the LHS vanishes for $i<M$. Therefore the only map of interest is
\begin{equation}\label{zh0}
H^M(P,{D}(\underline{\ml}(-a))[-2M]) \to H^M(P,\underline{\ml}(a)).
\end{equation}

Therefore we have a canonical element
\begin{equation}\label{portation}
\Sigma\in H^M(P,\underline{\ml}(-a))\otimes H^M(P,\underline{\ml}(a))= H^{2M}(P\times P, \underline{\ml}(a)\boxtimes \underline{\ml}(-a)).
\end{equation}

Now $A^{2M}(U\times U)=A^M(U)\tensor A^M(U)$ surjects onto $H^{2M}(P\times P, \underline{\ml}(a)\boxtimes \underline{\ml}(-a))$, and we will show that $\Sigma$ is represented by an explicit form as follows:

\begin{theorem}\label{main}
Let (as in \eqref{diagonalSIntro}, with $f^{(1)}_i(x,y)=f_i(x)$ and $f^{(2)}_i(x,y)=f_i(y)$)
\begin{equation}\label{diagonalS}
S=\sum_{1\leq i_1<\dots <i_{M}\leq r}\prod_{s=1}^{M} a_{i_s}\dlog f^{(1)}_{i_s}\dlog f^{(2)}_{i_s}\in A^{2M}(U\times U).
\end{equation}
Then, there exists a non-zero constant $c\in\Bbb{C}$ such that the image of $cS$ in  $H^{2M}(P\times P, \underline{\ml}(a)\boxtimes \underline{\ml}(-a))$ equals $\Sigma$.
\end{theorem}
In the case of small weights, Theorem \ref{main} specializes to Theorem \ref{maino}.
\begin{remark}\label{saturdayy}
We can, for arbitrary weights, also  consider the map $H^M_c(U,\ml(a))\to H^M(U,\ml(a))$ which is easily see to factor through
\eqref{zh0}. The form $S$ induces an element in $H^{2M}(U\times U,\ml(a)\boxtimes\ml(-a))$, which by Theorem \ref{main},  represents $H^M_c(U,\ml(a))\to H^M(U,\ml(a))$. Therefore, Theorem \ref{maino} holds for arbitrary weights.
\end{remark}

\subsection{A restatement of Theorem 48}

Consider the diagram
\begin{equation}\label{zhGCC}
\xymatrix{
H^M(P,{D}(\underline{\ml}(-a))[-2M])\ar[dr]^{\alpha} \ar[r]^{\sim} & H^{M}(P, \underline{\ml}(-a))^*\ar[d]\ar[r]^{\sim} & \bigl(\frac{A^M(U)}{\eta\wedge A^{M-1}(U))}\bigr)^*\ar[d]^{S}\\
& H^M(P,\underline{\ml}(a)) & \frac{A^M(U)}{\eta\wedge A^{M-1}(U)}. \ar[l]^{\sim}
}
\end{equation}
Here $\alpha$ is the map \eqref{zh}. The triangle on the left commutes because of the commutative diagram \eqref{diagramo}. Theorem \ref{main} implies that the square on the right commutes. The vertical map on the right is the composite:
\begin{equation}\label{positecomp}
\bigl(\frac{A^M(U)}{\eta\wedge A^{M-1}(U))}\bigr)^*\to A^M(U)^*\leto{S} A^M(U)\to \bigl(\frac{A^M(U)}{\eta\wedge A^{M-1}(U))}\bigr).
\end{equation}

\subsection{Action of symmetries}
Suppose a finite group $G$ acts on the hyperplane arrangement, permuting the hyperplanes, and preserving the corresponding weights. Then $G$ acts on $P\times P$ as well, and preserves the degree $2M$ form $S$ (since any two forms commute) in the statement of Theorem \ref{main}. All objects and maps in \eqref{zhGCC} are preserved under the action of this finite group $G$.

\subsection{Change of compactification}\label{changeC}
Let $P'$ be another compactification of $P$, and we assume as before that there is a regular birational morphism $\pi:P'\to P$. We use notation from Section \ref{officehours} in this section. There is a natural commutative diagram, note $D$ commutes with $R\pi_*$, and maps \eqref{canon} and \eqref{canonRR}. The commutativity is because of adjunction properties again.
\begin{equation}
\xymatrix{
{D}(\underline{\ml}(-a))[-2M])\ar[r]\ar[d] & \underline{\ml}(a)\\
R\pi_*(D(\underline{\ml'}(-a)))\ar[r] & R\pi_*(\underline{\ml'}(a)).\ar[u]
}
\end{equation}
Therefore the image of the map \eqref{zhGC} does not depend upon the compactification chosen, under the identification of  $H^{\starr}(P,\underline{\ml}(a))$ with the cohomology of the Aomoto complex $(A^{\starr}(U),\eta\wedge)$.
The following is now immediate,
\begin{lemma}\label{november1}
The image of the map \eqref{zhGC} inside the cohomology of the Aomoto complex is
$$\{x\in \frac{A^M(U)}{\eta\wedge A^{M-1}(U)}\mid \exists\  \tau\in A^{M}(U)^*,\  \tau(\eta\wedge A^{M-1}(U))=0,\ x=[S(\tau)]\}.$$
Therefore the image, under this identification, does not change if each $a_i$ is multiplied by multiplied by the same non-zero scalar.
\end{lemma}
\subsubsection{Representation of image by algebraic forms}\label{slika}
\begin{question}
Is the image of \eqref{zhGC} for $i=M$ (the image is zero otherwise), the linear span of $[\Omega]$ with $\Omega\in A^M(U)$ such that $\Omega$ does not have poles on divisors $E_{\alpha}$ with $a_{\alpha}=0$?
 \end{question}
\begin{lemma}\label{linearqQ}
The image of  \eqref{zhGC} for $i=M$  contains the linear span of $[\Omega]$ with $\Omega\in A^M(U)$ such that $\Omega$ does not have poles on divisors $E_{\alpha}$ with $a_{\alpha}=0$.
\end{lemma}
 \begin{proof}
 Let $\Omega$ be as in the statement of the lemma. By Lemma \ref{november1}, we are allowed to replace the numbers $a_i$ by $ca_i$ where $|c|$ is  very small. We then want to show the following:
 \begin{itemize}
 \item The element $[\Omega]\in H^M(U,\ml(a))$ is in the image of the mapping $H^M_c(U,\ml(a))\to H^M(U,\ml(a))$.
 \end{itemize}

 Let $V=P-\cup E'_{\alpha}$ where the union is over $\alpha$ such that $a_{\alpha}\neq 0$. Let $j:U\to V$. Then $H^M_c(U,\ml(a))= H^M(V,j_!{\ml(a)})$ and $H^M(U,\ml(a))=H^M(V, Rj_*\ml(a))$. By \cite{EV}, and Remark \ref{todayy} there is a factorization
  $$H^M(V,j_!{\ml(a)})\leto{\sim} H^M(V,\Omega^{\starr}(\log E)(-E))\to H^M(V, \Omega^{\starr}(\log E))\leto{\sim} H^M(V, Rj_*\ml(a)).$$

Now any $\Omega$ with the pole property above gives a class in  $H^0(V,\Omega^M(\log E)(-E))$ which is $d+\eta$ closed, and hence an element in $H^M(V,\Omega^{\starr}(\log E)(-E))$. The corresponding element of $H^M(V, \Omega^{\starr}(\log E))=H^M(V, Rj_*\ml(a))$ is the image of $[\Omega]$ under the isomorphism
$\frac{A^M(U)}{\eta\wedge A^{M-1}(U)}\to Rj_*\ml(a)$.
\end{proof}
In the case of arrangements coming from representation theory, the reverse implication in the above question is also true (for suitable isotypical components), see Section \ref{labelo} and Proposition \ref{reverso}.





\section{Steps in the proof of Theorem 48}
\subsection{The product hyperplane arrangement}
On $W\times W$ define polynomials of degree one ${f}^{(1)}_i(x,y)=f_i(x)$ and ${f}^{(2)}_i(x,y) =f_i(y)$. This defines a hyperplane arrangement $\mathcal{C}'$ in $W\times W$ given by hyperplanes ${f}^{(b)}_i=0, b=1,2,\ i=1,\dots,r$. We get a weighted
arrangement by assigning the weights $\tilde{a}$: assign $a_i$ to ${f}^{(1)}_i(x,y)=0$ and $-a_i$ to  ${f}^{(2)}_i(x,y) =0$.
The associated differential form is
\begin{equation}\label{SEl}
{\eta'}=\sum_{i=1}^r a_i (\dlog{f}^{(1)}_i -\dlog{f}^{(2)}_i)=\eta^{(1)}-\eta^{(2)},
\end{equation}
 where $\eta^{(b)}=\sum_{i=1}^r a_i \dlog{f}^{(b)}_i,\ b=1,2.$

 Note that $P\times P-U\times U$ is a divisor with normal crossings $E'$. The resulting object in the derived category of $P\times P$, i.e., $\underline{\ml}(\tilde{a})$ equals $\underline{\ml}(a)\boxtimes \underline{\ml}(-a)$, this can be seen also by showing the corresponding Aomoto cohomology for $P\times P$ equals the external product of the Aomoto complexes.

There is a switching factors automorphism $\sigma$ of $P\times P$
which preserves the arrangements, but is $-1$ on weights of hyperplanes: $\sigma(x,y)=(y,x)$.

\subsection{Main steps}\label{outline}

Let $\mk=\underline{\ml}(a)\boxtimes \underline{\ml}(-a)\in D^b_c(P\times P)$.
\begin{enumerate}
\item We will construct a non-zero class
$$\delta_{\alpha}\in H^{2M}_{\Delta}(P\times P, \mk)$$
which maps to $\Sigma\in H^{2M}(P\times P, \mk)$. The details can be found in Section \ref{breakfast}. Recall that $\Sigma$ was defined in Section \ref{afternoon} (see \eqref{portation}). We also explain in Section \ref{breakfast} that  it suffices to show that the image of $\delta_{\alpha}$ in $H^{2M}(P\times P, \mk)$ coincides with the class $[S]$.\\
\item  In Lemma \ref{tif}, we show that $H^{2M}_{\Delta}(P\times P, \mk)=\Bbb{C}$, and that $\delta_{\alpha}$ is a generator of this group.\\
    \item We construct a {\em non-zero} element $S_{\Delta}\in H^{2M}_{\Delta}(P\times P, \mk)$ which maps to $[S]\in H^{2M}(P\times P, \mk)$, where $S$ is as in Theorem \ref{main}. We  show that $S_{\Delta}$ is a non-zero multiple of $\delta_{\alpha}$. We refer to Section \ref{hepburn} for details.
\end{enumerate}
The following consistency check is a good way of viewing the third step above. For such a $S_{\Delta}$ to exist,  the image of
$[S]$ in $H^{2M}(P\times P-\Delta, \mk)$ needs to be zero since we have an exact sequence (see the exact sequence on cohomology induced by the distinguished exact triangle in \cite[Exercise I.2.5]{KS} with $S=S_{k-1}$ and $S_k=\emptyset$)
$$H^{2M}_{\Delta}(P\times P, \mk)\to H^{2M}(P\times P,\mk)\to  H^{2M}(P\times P-\Delta,\mk).$$

For a linear function $F$ on $W$, we can form $U_F=P\times P-\overline{{Z}(h)}\subseteq P\times P-\Delta$ where $h(x,y)= F(x)-F(y)$, and
$Z(h)\subset U\times U$ is the zero set of $h$. Therefore
the image of $[S]$ in $H^{2M}(U_F,\mk)$ should be zero. This would be true if the image of $$S\in A^{2M}\bigl(U_F\cap (U\times U)\bigr)\subseteq H^0(U_F,\Omega^{2M}_{P\times P}(\log E'))$$ is exact, note that $U_F-E'= U_F\cap (U\times U)$. This follows from Corollary \ref{grundlegendSpeciale}, which implies (see \eqref{diagonalgeneral} for the definition of $S^{(M-1)}$)
\begin{equation}\label{seventyy}
S = (\eta^{(1)}-\eta^{(2)})\wedge \dlog h \wedge S^{(M-1)}.
\end{equation}
Varying $F$, one could hope that $U_F$ form an open cover of
$P\times P-\Delta$, and that a generalization of \eqref{seventyy}  would show that $[S]$ vanishes in $H^{2M}(P\times P-\Delta, \mk)$.  The element $S_{\Delta}$ can then be hoped to arise from a cone construction.

But it seems to be difficult to ensure that $U_F$ form an open cover of
$P\times P-\Delta$. Instead we use a slightly different strategy: We suitably blow up $P\times P$ outside of $U\times U$ and use a similar argument.

\section{Some generalities}\label{breakfast}
Let $\alpha:F\to G$ be a morphism in $D^b_c(P)$. We will use the considerations of this section with $\alpha$ the mapping \eqref{zh}. First note that $\alpha$ produces a cohomology class
\begin{equation}\label{deltaforce}
\delta_{\alpha}\in H^0\bigl(P,\Gamma_{\Delta} ({{{{\srhom}}}}(p_1^{-1}F, p_2^{!}G))\bigr)
\end{equation}
by Proposition 3.1.14 of [KS]. Clearly $\delta_{\alpha}$ also produces an element in  $H^0$ of ${\srhom}(p_1^{-1}F, p_2^{!}G)$, which is the object on the top right of the following diagram (the remaining objects and maps are explained below):
\begin{equation}\label{diagramo}
\xymatrix{R\Gamma(P\times P,DF\boxtimes G) \ar[r] & {\operatorname{Rhom}}(p_1^{-1}F, p_2^{!}G)\ar[d]\\
R\Gamma(P,DF)\tensor R\Gamma(G)\ar[u] & {{\operatorname{Rhom}}}((Ra_P)_*F,(Ra_P)_*G)\\
R\Gamma(P,F)^*\tensor R\Gamma(G).\ar[u]\ar[ur]
}
\end{equation}
 In the rest of this section we recall the maps in \eqref{diagramo} which are all isomorphisms. The vertical map on the right in \eqref{diagramo} is the  morphism as in Proposition 3.1.15 of [KS]:
\begin{lemma}There is a natural isomorphism
\begin{equation}\label{promo}
 {{\operatorname{Rhom}}}(p_1^{-1}F, p_2^{!}G)\leto{\sim} {\operatorname{RHom}}((Ra_P)_!F,(Ra_P)_*G).
 \end{equation}
\end{lemma}

The map \eqref{promo} arises as follows
 $$ {{\operatorname{Rhom}}}(p_1^{-1}F, p_2^{!}G)=(Ra_P)_* (Rp_2)_*{{{{\srhom}}}}(p_1^{-1}F, p_2^{!}G)\to (Ra_P)_*{{{{\srhom}}}}((Rp_2)_{!}p_1^{-1}F, G),$$
which in  turn maps to
$$(Ra_P)_*{{ {{\srhom}}}}(a_P^{-1}(Ra_P)_!F,G)\to {{\operatorname{Rhom}}}((Ra_P)_!F,(Ra_P)_*G).$$

We note that the image of $\delta_{\alpha}$ (defined in \eqref{deltaforce}) under \eqref{promo} is the map induced by $\alpha$:
$$(Ra_P)_{!}F \leto{\alpha} (Ra_P)_{!}G\to (Ra_P)_* G.$$

\begin{lemma}\label{derived}
\begin{equation}\label{muse}
DF\boxtimes G \leto{\sim}{{ {{\srhom}}}}(p_1^{-1}F, p_2^{!}G).
\end{equation}
\end{lemma}

 We recall  how the morphism in   Lemma \ref{derived} is constructed
$$DF\boxtimes G\to p_1^{-1}{{ {{\srhom}}}}(F, a_P^{!}\Bbb{C})\otimes p_2^{-1}G.$$
Using equation (2.6.27), page 114 of [KS],
\begin{equation}\label{compose}
p_1^{-1}{{{{\srhom}}}}(F, a_P^{!}\Bbb{C})\otimes p_2^{-1}G\to
{{{{\srhom}}}}(p_1^{-1}F, p_1^{-1}a_P^{!}\Bbb{C})\otimes p_2^{-1}G\to {{{{\srhom}}}}(p_1^{-1}F, p_1^{-1}a_P^{!}\Bbb{C}\otimes p_2^{-1}G).
\end{equation}

Using prop 3.1.9, (iii), page 146 of [KS] we have a map induced by adjunction:
$$p_1^{-1}a_P^{!}\Bbb{C}\to p_2^{!} a_P^{-1}(\Bbb{C})\to p_2^{!}\Bbb{C}.$$

By Proposition 3.1.11, on page 147 of [KS], there is a natural map $p_2^{!}\Bbb{C}\otimes p_2^{-1}G\to p_2^{!}G$. We compose \eqref{compose} with this morphism to complete the argument. Proposition 3.4.4  in [KS] shows that the constructed map is an isomorphism.

\begin{lemma}\label{beanie}
Let $F'=DF,G\in D^b_c(P)$. Then
$$ R\Gamma(P,F')\tensor R\Gamma(G)\leto{\sim} R\Gamma(P\times P,F'\boxtimes G)$$
and
$R\Gamma(P,F')= R\Gamma(P,T)^*.$
\end{lemma}
The morphism in Lemma \ref{beanie} arises as follows:
$$(Ra_P)_*((Rp_2)_*p_1^{-1}F'\otimes G)\leto{\sim}(Ra_P)_*(Rp_2)_*(p_1^{-1}F'\otimes p_2^{-1} G)\leto{\sim} R\Gamma (F'\boxtimes G).$$

Now $ Ra_P^*(Ra_P)_* F'\to (Rp_2)_*p_1^{-1}F'$. Therefore we have an isomorphism map
$$(Ra_P)_*(Ra_P^*(Ra_P)_* F'\otimes G)\to R\Gamma(P\times P,F'\boxtimes G)$$
and hence $$(Ra_P)_*F'\otimes (Ra_P)_* G\to R\Gamma(F'\boxtimes G).$$

Therefore all maps in \eqref{diagramo} have been constructed, the key claim (standard) is
\begin{lemma}
The diagram \eqref{diagramo} commutes.
\end{lemma}

Now let $\alpha$ be the mapping \eqref{zh}.
We therefore have a class
$$\delta_{\alpha}\in H^{2M}_{\Delta}(P\times P, \underline{\ml}(a)\boxtimes \underline{\ml}(-a))$$
which maps to $\Sigma\in H^{2M}(P\times P, \underline{\ml}(a)\boxtimes \underline{\ml}(-a))$, which in turn induces the mapping \eqref{zhGC}.

\begin{lemma}\label{tif}
\begin{enumerate}
\item[(a)]  Setting $\tilde{k}:\Delta\cap (U\times U)\to {\Delta}$.$$i_{{\Delta}}^{!}  (\underline{\ml}(a)\boxtimes\underline{{\ml}}(a))= R\tilde{k}_{*}\Bbb{C}[2M].$$
\item[(b)] $H^{2M}_{\Delta}(P\times P, \underline{\ml}(a)\boxtimes \underline{\ml}(-a))=\Bbb{C}$.
\item[(c)] Under the isomorphism in (b), $\delta_{\alpha}$ corresponds to $1\in \Bbb{C}$.
\end{enumerate}
\end{lemma}
\begin{proof}
 We claim that for all points $p\in {\Delta}- U\times U$, there is a $\beta$ such that $p\in {E'}_{\beta}$ with ${a'}_{\beta}\in \Bbb{C}-\{1,2,\dots\}$.

  This is true because  $\sigma$ fixes $p$, and if $p\in {E'}_{\beta}$ then $p\in {E'}_{\sigma(\beta)}$. Therefore the claim follows from the identities,
   ${a'}_{\sigma(\beta)}=-{a'}_{\beta}$ ($\sigma(\eta')=-\eta'$, and $\sigma({E'}_{\beta})={E'}_{\sigma(\beta)}$ by definition of $\sigma(\beta)$).

 We can use standard  base change properties  (Proposition 3.1.9 on page 145, and Proposition 3.1.1 on page 147 of \cite{KS}). We also use the fact that ${\widetilde{\ml}}(a)$ restricted to $\Delta\cap (U\times U)$ is trivial. This implies (a), and (b) follows from (a). To prove (c) we may localize in $P$, and reduce to a standard compatibility. We only need that the factor is non-zero which is clear because otherwise the map $\alpha$ in question, i.e., the map \eqref{zh} would be zero on $U$.
\end{proof}

\section{Cohomology class of the diagonal}\label{classo}
Let $\mathcal{K}=\underline{\ml}(a)\boxtimes \underline{\ml}(-a)\in D^b_c(P\times P)$.
Our aim in this section (see Proposition \ref{diagonalelement}) is to construct an element $S_{\Delta}\in H^{2M}_{\Delta}(P\times P, \mathcal{K})$ which maps to
$[S]\in H^{2M}(P\times P, \mathcal{K})$ (see Theorem \ref{main} for the definition of $S$).

In this section we construct a suitable $\widetilde{P}$ birational to $P\times P$. Set $\widetilde{\Delta}\subset \wtp$ the strict transform of $\Delta$, and let $\underline{\widetilde{\ml}}(a)$ denote the object in the derived category of sheaves on $\widetilde{P}$ for the product arrangement (with weights $a$ and $-a$). By  Section \ref{officehours}, we know that $\underline{\widetilde{\ml}}(a)$ is represented by the Aomoto complex. We show
\begin{enumerate}
\item $[S]\in H^{2M}(\wtp,\underline{\widetilde{\ml}}(a))$ goes to zero when restricted to $\wtu=\wtp-\wtd$.
\item Carrying out the previous step keeping track of forms that appear in the vanishing, we construct an element in $H^{2M}_{\wtd}(\wtp,\underline{\widetilde{\ml}}(a))$.
\item  $H^{2M}_{\wtd}(\wtp,\underline{\widetilde{\ml}}(a))$ maps to $H^{2M}_{\Delta}(P\times P, \mathcal{K})$. We then define $S_{\Delta}$ to be the image of the element in the previous step.
\end{enumerate}

\subsection{Functorial Resolution and the diagonal}
Let $F_1=0,\dots,F_{M}=0$ be linear hyperplanes in $W$ such that $dF_1,\dots, dF_M$ are linearly independent.  The diagonal in $W\times W$ is cut out by $h_j(x,y)=F_j(x)-F_j(y)=0$ with $j=1,\dots,M$. These meet transversally in $W\times W$. Let $Z(h_j)$ be the variety of zeroes of $h_j$ on $U\times U$, and $\overline{{Z}(h_j)})\subset P\times P$ its closure.
\begin{remark}
 The intersection $\cap_{j=1}^M\overline{{Z}(h_j)}\subset P\times P$ may be bigger than the diagonal $\Delta\subset P\times P$, and therefore the complements
$P\times P- \overline{{Z}(h_j)}$ may not cover $P\times P-\Delta$. We pass to a blowup of $P\times P$ such that the strict transforms of $\overline{{Z}(h_j)}$ meet properly. It then follows that the strict transform meet in a locus which maps to $\Delta\subset P\times P$.
\end{remark}
Let $\pi:\widetilde{P}\to P\times P$ be a functorial resolution of singularities \cite{BVP}, following Section \ref{TonyS}, of the pair $(P\times P, E'\cup \bigcup_{j=1}^M \overline{{Z}(h_j)})$, where $E'=P\times P - U\times U$. Since the zero loci $Z(h_j)$ meet transversally in $U\times U$, by property (d) in Section \ref{TonyS},
$\pi$ is an isomorphism over $U\times U$. Set
$$\widetilde{Z}(h_j)=\overline{\{(x,y)\in U\times U: F_j(x)= F_j(y)\}}\subseteq \widetilde{P}.$$
Let $\widetilde{E}=\pi^{-1}(E'\cup\bigcup_{j=1}^M \overline{{Z}(h_j)})$, a divisor with simple normal crossings. Write $\widetilde{E}$ as a union $\cup_{\beta} \widetilde{E}_{\beta}$, let $\tilde{a}_{\beta}$ be the residue of $\eta'=\eta^{(1)}-\eta^{(2)}$ along $\widetilde{E}_{\beta}$.
As will be clear from what follows, $\widetilde{P}$ is a much better place for actual computations.
\begin{defi}\label{scary}
It is clear that $\widetilde{Z}(h_j)$ are irreducible components of $\widetilde{E}$. Write
$$\widetilde{E}=\widetilde{E}_0\cup (\bigcup_{j=1}^M \widetilde{Z}(h_j))$$
where $\widetilde{E}_0$ is the union of the other irreducible components.
\end{defi}
\begin{itemize}
\item The automorphism  $\sigma$ of switching the two factors in $P\times P$ lifts to $\widetilde{P}$, so that each $\widetilde{E}_{\beta}$ goes to $\widetilde{E}_{\beta'}$ with $\tilde{a}_{\beta'}+\tilde{a}_{\beta}=0$
    (canonical resolution of singularities).
\item Let $\widetilde{\Delta}$ be the closure of $\Delta\cap (U\times U)$ in $\widetilde{P}$. It is easy to see that $\cap_{j=1}^M \widetilde{Z}(h_j)  = \widetilde{\Delta}$ using transversality (and dimension counting). $\widetilde{\Delta}$ maps to $\Delta$ under $\widetilde{P}\to P\times P$.
\end{itemize}
Let $\underline{\widetilde{\ml}}(a)$ denote the object in the derived category of sheaves on $\widetilde{P}$ for the arrangement, (without the hyperplanes $h_j=0$), and the form $\sum_{i=1}^r a_i (\dlog{f}^{(1)}_i(x,y)-\dlog{f}^{(2)}_i(x,y))$, this is quasi-isomorphic to  $(\Omega^{\starr}_{{\wtp}}(\log \widetilde{E}_0),d+\eta^{(1)}-\eta^{(2)})$. Recall the notation Let $\mathcal{K}=\underline{\ml}(a)\boxtimes \underline{\ml}(-a)$.

By Section \ref{canon}, we have a quasi-isomorphism
\begin{equation}\label{YesM}
R\pi_*\underline{\widetilde{\ml}}(a)\to \mk.
\end{equation}
\begin{equation}\label{diagramo2}
\xymatrix{
H^{2M}_{\widetilde{\Delta}}(\widetilde{P}, \underline{\widetilde{\ml}}(a))\ar[r]\ar[d] & H^{2M}(\widetilde{P}, \underline{\widetilde{\ml}}(a))\ar[r]\ar[d] & H^{2M}(\wtu, \underline{\widetilde{\ml}}(a))\ar[d]\\
 H^{2M}_{\Delta}(P\times P, \mathcal{K})\ar[r] &  H^{2M}(P\times P, \mathcal{K})\ar[r] & H^{2M}(P\times P-\Delta, \mathcal{K}).
}
\end{equation}

Recall that $\wtu=\wtp-\wtd$.

The middle vertical arrow is an isomorphism. The vertical arrow on the far right is obtained by restricting \eqref{YesM} to $P\times P-\Delta$, using the restriction map
$$H^{2M}(\wtu, \underline{\widetilde{\ml}}(a))\to H^{2M}(\pi^{-1}(P\times P-\Delta), \underline{\widetilde{\ml}}(a)).$$
The vertical map on the far left is induced by the isomorphism \cite[Proposition 3.1.9]{KS}
 $$R\pi'_*i_{\pi^{-1}\Delta}^{!}\underline{\widetilde{\ml}}(a)\to i_{\Delta}^{!}R\pi_*\underline{\widetilde{\ml}}(a)$$
and the inclusion $\widetilde{\Delta}\subseteq \pi^{-1}\Delta$.  The following lemma, which uses property (e) of resolution of singularities from Section \ref{TonyS}, shows that this vertical map is an isomorphism. This result is not used, but gives
a more satisfying picture.
\begin{lemma}
$i_{\widetilde{\Delta}}^{!} \underline{\widetilde{\ml}}(a))= Rk'_{*}\Bbb{C}[2M]$ where $k':\Delta\cap (U\times U)\to \widetilde{\Delta}$.
\end{lemma}
\begin{proof}
 We claim that for all points $p\in \widetilde{\Delta}- U\times U$, there is a $\beta$ such that $p\in \widetilde{E}_{\beta}$ with $\tilde{a}_{\beta}\in \Bbb{C}-\{1,2,\dots\}$.

  This is true because $\sigma:\widetilde{P}\to\widetilde{P}$ fixes $p$, and if $p\in \widetilde{E}_{\beta}$ then $p\in \widetilde{E}_{\sigma(\beta)}$. Therefore the claim follows from the identity
   $$\tilde{a}_{\sigma(\beta)}=- \tilde{a}_{\beta}$$
   ($\sigma(\eta')=-\eta'$, and
   $\sigma(\widetilde{E}_{\beta})=\widetilde{E}_{\sigma(\beta)}$ by definition of $\sigma(\beta)$).

 We have used standard  base change properties  (Proposition 3.1.9 on page 145, and Proposition 3.1.1 on page 147 of \cite{KS}), also that $\underline{\widetilde{\ml}}(a)$ restricted to $\Delta\cap (U\times U)$ is trivial.
\end{proof}
\begin{proposition}\label{varnish}
$S\in H^M(\wtp,\underline{\widetilde{\ml}}(a))$ goes to zero in $H^M(\wtu,\underline{\widetilde{\ml}}(a))$, where $\wtu=\wtp-\wtd$.
\end{proposition}

Now, $\wtu=\widetilde{P}-\widetilde{\Delta}$ has an open covering by open subsets of the form
$\widetilde{P}-\widetilde{Z}(h_j).$ An intersection of these open subsets has the form
$$\widetilde{P}_J=\widetilde{P}-\bigcup_{j\in J}\widetilde{Z}(h_j)=\bigcap_{j\in J} \bigl(\widetilde{P}-\widetilde{Z}(h_j)\bigr).$$
where $J\subset\{1,\dots,M\}$.

\subsection{\Ccech complexes}\label{check}
For a sheaf $\mathcal{F}$ on $\wtu=\widetilde{P}-\wD$, let
$$C^p(\mf)=\prod_{j_0<\dots<j_p} \mf(\widetilde{P}_{\{j_0,\dots,j_p\}}).$$

Elements $\alpha\in C^p(\mf)$ are determined by
giving  elements
$$\alpha_{j_0,\dots,j_p}\in  \mf(\widetilde{P}_{\{j_0,\dots,j_p\}}).$$

Define $\delta:C^p\to C^{p+1}$ by
$$(\delta\alpha)_ {j_0,\dots,j_{p+1}}= \sum_{k=0}^{p+1} (-1)^k \alpha_{j_0,\dots,\hat{j}_k,\dots j_{p+1}}.$$

Also form the sheaf version (with a corresponding differential)
$$\mathcal{C}^p(\mf)=\prod_{j_0<\dots<j_p} (k_J)_*\mathcal{F}\mid \widetilde{P}_{\{j_0,\dots,j_p\}},$$
where $k_J$ denotes the inclusion of  $\widetilde{P}_{\{j_0,\dots,j_p\}}\to \widetilde{P}-\wD$.
It is known that there is a quasi-isomorphism $\mf\to (\mathcal{C}^p(\mf),\delta)$ (see \cite[Lemma III.4.2]{Hartshorne}).

Now suppose $(\mf^{\starr},d_{\mf})$ is a complex of sheaves. Let $(\mathcal{C}^{\starr}(\mf^{\starr},d_{\mf}),D)$ be the complex with
$$\mathcal{C}^n(\mf^{\starr},d_{\mf})=\oplus_{p+q=n} \mathcal{C}^p(\mf^q)$$
and differential given by $D=d_{\mf}+\delta$.
It is again known that there is a quasi-isomorphism $\mf^{\starr}\to \mathcal{C}^{\starr}(\mf^{\starr})$.

\subsection{\Ccech complex for twisted log de Rham complex }\label{defD}
Let, as before,  $\widetilde{U}=\widetilde{P}-\widetilde{\Delta}$. Now,  $\widetilde{U}$ has a covering   $\widetilde{P}-\widetilde{Z}(h_j)$.
Let  $\mf^{\starr}=(\Omega^{\starr}_{\widetilde{P}}(\log \widetilde{E}_0),d+\eta^{(1)}-\eta^{(2)})$ which is quasi-isomorphic to $\underline{\widetilde{\ml}}(a)$.
($\widetilde{E}_0$ was defined in Definition \ref{scary}).
Let $\mf^{\starr}\to \mathcal{I}^{\starr}$ be a quasi-isomorphism where $\mathcal{I}^{\starr}$ is a complex of injectives. We obtain a commutative diagram of complexes
$$\xymatrix{
  &                               \Gamma(\wtu,\mathcal{C}^{\starr}(\mf^{\starr}))\\
\Gamma(\wtp,\mf^{\starr})\ar[r]\ar[ur] & \Gamma(\wtu,\mf^{\starr}).\ar[u]
}
$$
Here $\mathcal{C}^{\starr}(\mf^{\starr})$ is the \Ccech complex of $\mf^{\starr}$ on $\wtu$, see Section \ref{check} for the definition.
The above diagram maps to a similar diagram associated to $\mathcal{I}^{\starr}$: The vertical map in the diagram below is a quasi-isomorphism, and $\mathcal{C}^{\starr}(\mI^{\starr}))$ is a complex of injectives.
$$\xymatrix{
  &                               \Gamma(\wtu,\mathcal{C}^{\starr}(\mI^{\starr}))\\
\Gamma(\wtp,\mI^{\starr})\ar[r]\ar[ur] & \Gamma(\wtu,\mI^{\starr}).\ar[u]
}
$$
Hence we obtain
\begin{equation}\label{newman}
\xymatrix{
\Gamma(\wtp,\mf^{\starr})\ar[r]\ar[d] & \Gamma(\wtu,\mathcal{C}^{\starr}(\mf^{\starr}))\ar[d]\\
\Gamma(\wtp,\mI^{\starr})\ar[r] & \Gamma(\wtu,\mathcal{C}^{\starr}(\mI^{\starr})).
}
\end{equation}

We are given $[S]\in H^{2M}(\Gamma(\wtp,\mf^{\starr}))=H^{2M}(\Gamma(\wtp,\mI^{\starr}))$ . Therefore, to prove Proposition \ref{varnish}, it suffices to show that the image of $[S]$ in $H^{2M}(\Gamma(\wtu,\mathcal{C}^{\starr}(\mf^{\starr})))$ is zero. Note that
$H^{2M}\Gamma(\wtu,\mathcal{C}^{\starr}(\mI^{\starr}))=H^{2M}(\wtu,\mf^{\starr})$.

Let $$D^n =\oplus_{p+q=n} \oplus _J A^q(U_J)$$
where $J$ runs through $J=\{j_0<\dots<j_p\}$
and $$U_J =\wtp_J\cap (U\times U)=\{(x,y)\mid \forall i\in J, F_i(x)\neq F_i(y)\}.$$

Note that $D^n\subset \mathcal{C}^n(\mf^{\starr})$: This is because for any $J$,
$$A^q(U_J)\subset H^{0,\op{an}}(\wtp_J,\Omega^q_{\wtp}(\widetilde{E}_0)).$$
The last inclusion needs an explanation: $\wtp$ is a compactification of $U_J$, with complement $$D_J=\widetilde{E}_0\cup \bigcup_{j\in J} \widetilde{Z}(h_j). $$
Therefore, $A^q(U_J)=H^0(\wtp,\Omega^q_{\wtp}(D_J))$, justifying the inclusion above.

\subsection{Some log forms}
The following Shapovolov form, a log form on $U\times U$ was defined in equation \eqref{diagonalSIntro}
\begin{equation}\label{diagonal}
S=\sum_{1\leq i_1<\dots <i_{M}\leq r}\prod_{s=1}^{M} a_{i_s}\dlog f^{(1)}_{i_s}\dlog f^{(2)}_{i_s}.
\end{equation}
In fact for any $0\leq b\leq M$, we may define log forms on $U\times U$ of degree $b$:
\begin{equation}\label{diagonalgeneral}
S^{(b)}=\sum_{1\leq i_1<\dots <i_{b}\leq r}\prod_{s=1}^{b} a_{i_s}\dlog f^{(1)}_{i_s}\dlog f^{(2)}_{i_s}
\end{equation}
so that $S^{(0)}=1$ and $S^{(M)}=S$.

Further for every $\{q_1<\dots<q_w\}\subseteq [1,M]$, and $b$ such that $w+b=M$, we define a log form on $U\times U -\cup_{j=1}^w Z(F^{(1)}_{q_j}-F^{(2)}_{q_j})$. Here $F^{(1)}_j(x,y)=F_j(x)$ and $F^{(2)}_j(x,y)=F_j(y)$, note that we had earlier
defined $h_j(x,y)=F_j(x)-F_j(y)= F^{(1)}_j(x,y)- F^{(2)}_j(x,y)$.

$$S_{q_1\dots q_w}:=S^{(b)}\prod_{j=1}^w \dlog (F^{(1)}_{q_j}- F^{(2)}_{q_j})$$
a form of degree $2b+w= 2(M-w)+w=2M-w$. When $w=0$, we recover $S$.
\begin{remark}\label{addington}
The form $S_{1,2,\dots,M}$ corresponds to $w=M$  and $b=0$, hence
$$S_{1,2,\dots,M}:=\prod_{j=1}^M \dlog (F^{(1)}_{j}- F^{(2)}_{j}).$$
This form does not depend upon the weights $a$.
\end{remark}
\subsection{Proof of Proposition 74}\label{varnishproof}

The Shapovolov element $S|U_{\{i\}}$, and zero in all intersections gives a closed element in $D^{2M}$.
 We claim that this is zero in cohomology. For this we need elements  $\alpha=\{\alpha_J\}$ in
$D^{2M-1} =\oplus_{s+q=2M-1} \oplus _J A^q(U_J)$ with $J=\{j_0<\dots<j_s\}$ where $U_J =\wP^J\cap (U\times U)$ (defined in Section \ref{defD}).
Define
\begin{equation}\label{outsourcinghomework}
\alpha_J=S_{j_0,j_1\dots, j_s}
\end{equation}
which is a degree $q=2M-(s+1)$ form. In particular when $J=\{1<2<\dots<M\}$, we get $\alpha_J=\prod_{i=1}^M\dlog (F^{(1)}_j-F^{(2)}_j)$.
Proposition \ref{grundlegend} below shows that the differential of $\alpha$ is the image of $[S]$ in $H^M(\Gamma(\wtu,\mathcal{C}^{\starr}(\mf^{\starr})))$, proving Proposition \ref{varnish}.
\begin{proposition}\label{grundlegend}
Let $1\leq k\leq M$, then
\begin{equation}\label{formule1}
\sum_{j=1}^{k}(-1)^{j+1} S_{q_1,\dots \hat{q}_j,\dots q_{k}}=(\eta^{(1)}-\eta^{(2)})\wedge S_{q_1\dots q_{k}}.
\end{equation}
\end{proposition}
We give a proof of Proposition \ref{grundlegend} in Section \ref{keyF}. We note the following corollary.

\begin{corollary}\label{grundlegendSpeciale}
For any $q=1,\dots,M$, $S = (\eta^{(1)}-\eta^{(2)})\wedge S_{q}$.
\end{corollary}
\subsection{Construction of the class supported on the diagonal from $S$}\label{dirty}
 Our aim in this section is to construct
$S'\in H^M_{\wtd}(\wtp,\underline{\widetilde{\ml}}(a))$ which maps to $S\in H^M(\wtp,\underline{\widetilde{\ml}}(a))=H^M(P\times P,\mk)$. The element $S_{\Delta}\in  H^M_{\Delta}(P\times P,\mk)$  will then be the image of $S'$, see diagram \eqref{diagramo2}.

Using notation from Section \ref{defD}, consider the exact sequences:
$$0\to \Gamma_{\wtd}(\wtp,\mathcal{I}^{\starr})\to \Gamma(\wtp,\mathcal{I}^{\starr})\leto{f} \Gamma(\wtu,\mathcal{I}^{\starr})\to 0.$$
The hypercohomology groups $H^{\starr}_{\wtd}(\wtp,\mathcal{F}^{\starr})$ are therefore canonically isomorphic to cohomology of the cone of $f$, which in turn is canonically isomorphic to the cohomology of the cone of $\Gamma(\wtp,\mathcal{I}^{\starr})\to \Gamma(\wtu,\mathcal{C}^{\starr}(\mathcal{I}^{\starr}))$. We therefore obtain, using the diagram \eqref{newman},
\begin{lemma}
The cohomology of the cone of
$\Gamma(\wtp,\mathcal{F}^{\starr})\to \Gamma(\wtu,\mathcal{C}^{\starr}(\mathcal{F}^{\starr}))$ admits a mapping to $H^{\starr}_{\wtd}(\wtp,\mathcal{F}^{\starr})$ consistent with various exact sequences.
\end{lemma}
\begin{remark}
We have used the following standard fact about cones of morphisms: Suppose
$$0\to A^{\starr}\to B^{\starr}\leto{f} C^{\starr}\to 0$$
is a short exact sequence of complexes (exact in each degree), then $A^{\starr}$ (with a shift) canonically maps to the cone of $f$, and this map is a quasi-isomorphism (see e.g., \cite[Exercise 1.5.7]{weibel}).
\end{remark}

The elements $S$ and $\alpha$ from Section \ref{varnishproof} give an explicit element in the cohomology of the cone of $\Gamma(\wtp,\mathcal{F}^{\starr})\to \Gamma(\wtu,\mathcal{C}^{\starr}(\mathcal{F}^{\starr}))$. We therefore get an element in
in $H^{\starr}_{\wtd}(\wtp,\mathcal{F}^{\starr})$. Let $S_{\Delta}\in H^{2M}_{\Delta}(P\times P, \mathcal{K})$ be the image of this element under the vertical map on the left in the diagram \eqref{diagramo2}.  Therefore,
\begin{proposition}\label{diagonalelement}
There is an element $S_{\Delta}\in H^{2M}_{\Delta}(P\times P, \mathcal{K})$ which maps to
$[S]\in H^{2M}(P\times P, \mathcal{K})$. Here $\mathcal{K}=\underline{\ml}(a)\boxtimes \underline{\ml}(-a)\in D^b_c(P\times P)$.
\end{proposition}

\section{Proof of Theorem 48}\label{hepburn}
\subsection{Some standard facts}
Consider the natural map $H^{2M-1}(\Bbb{A}^M-\{0\},\Bbb{C})\to H^{2M}_{0}(\Bbb{A}^M,\Bbb{C})$ coming from the standard supports sequence.
We have a natural map $H^{2M}_{0}(\Bbb{A}^M,\Bbb{C})=H^{2M}(i_0^{!}\Bbb{C})=\Bbb{C}$ because
of the standard isomorphism $i_0^{!}\Bbb{C}=\Bbb{C}[-2M]$. We hence obtain a map $H^{2M-1}(\Bbb{A}^M-\{0\},\Bbb{C})\to \Bbb{C}$.

We also have a map $H^{2M-1}(\Bbb{A}^M-\{0\},\Bbb{C})\to \Bbb{C}$ given by integration on the $2M-1$ sphere. This coincides with the map (up to a sign)
with the map in the previous paragraph.

There is a map $H^0(\Bbb{A}^M -\cup_{i=1}^M D_i, \Omega^M)\to H^{2M-1}(\Bbb{A}^M-\{0\},\Bbb{C})$, coming from the \Ccech covering of $\Bbb{A}^M-\{0\}$ by
$\Bbb{A}^M- D_i$ where $D_i=\{(z_1,\dots,z_M)\mid z_i=0\}$, and the holomorphic de Rham resolution of $\Bbb{C}$. The composite map $H^0(\Bbb{A}^M -\cup_{i=1}^M D_i, \Omega^M)\to H^{2M-1}(\Bbb{A}^M-\{0\},\Bbb{C})\to\Bbb{C}$
is $\frac{1}{(2\pi \sqrt{-1})^M}$ times the Grothendieck residue. We record the following consequence (see \cite[Page 651]{GH}):

\begin{lemma}\label{GHresult}
The topological map  $H^0(\Bbb{A}^M -\cup_{i=1}^M D_i, \Omega^M)\to \Bbb{C}$ sends $\frac{dz_1}{z_1}\dots\frac{dz_M}{z_M}$ to $(2\pi \sqrt{-1})^M$.
\end{lemma}
\subsection{Completion of the Proof of Theorem 48}
Recall from Section \ref{outline} the outline of the proof of Theorem \ref{main}. To complete the proof of Theorem \ref{main}, we need to show that $S_{\Delta}$ (constructed in Proposition \ref{diagonalelement}) is a non-zero multiple of the generator  $\delta_{\alpha}$ of $H^M_{\Delta}(P\times P,\underline{\ml}(a)\boxtimes\underline{\ml}(-a))$.

The construction of $S_{\Delta}$ can be localized, to any open subset of $U\times U\subset \widetilde{P}$ of the form $A\times A$ where $A$ is an  open polydisc contained in $U$. Therefore $S_{\Delta}$ gives an element in the $2M$-th cohomology cone $\op{Cone}_A(f)$ of
$$\Gamma(A\times A,\mathcal{F}^{\starr})\to \Gamma(A\times A-\Delta_A,\mathcal{C}^{\starr}(\mathcal{F}^{\starr})).$$

Consider the mapping $A\to A\times A$ given $\vec{z}\mapsto (\vec{z},0)$ where the forms $F_i$ correspond to the coordinates. This map pulls back the diagonal to $0\in A$. There is therefore a map from the $H^{2M}$ of $\op{Cone}_A(f)$ to $H^{2M}$ of the cone $\op{Cone}_{0}(f_0)$ of
$$\Gamma(A,\mathcal{F'}^{\starr})\to \Gamma(A-\{0\},\mathcal{C}^{\starr}(\mathcal{F'}^{\starr}))$$
where $\mathcal{F}'^{\starr}$ is the corresponding Aomoto complex of $A$. This map between the $H^{2M}$  of the cones is clearly an isomorphism, since the cones compute topological cohomology groups with supports.

Therefore we need to show that the forms in the construction of $S_{\Delta}$ when pulled back to
$A$ together give a non-zero  element in $H^{2M}_0(A,\ml(a))$. It is clear that $\dlog f_i^{(2)}$ pull back to zero, and hence the pull back element in the complex $\op{Cone}_{0}(f_0)$ sits in only one degree (see Remark \ref{addington} and Equation \eqref{outsourcinghomework}), and is $\prod_{i=1}^M \dlog z_i\in H^{M}(A-\cup_{i=1}^M \{z_i=0\},\Omega^M)$.
We need to multiply this function by a local generator of the local system $\ml(a)\times \ml(-a)_0$ and view it as an image of an element in $H^{2M-1}(A-\{0\},\Bbb{C})$ (computed in \Ccech cohomology of the holomorphic de Rham complex for the coordinate covering $A-\{z_i=0\}$ of $A-\{0\}$).

Let  $g$ is a holomorphic function on $A$ which is a local generator of $\ml(a)$, i.e., $dg+g\eta=0$.
The corresponding element of $H^{2M-1}(A-\{0\},\Bbb{C})$ is the image of $g g(0)^{-1}\prod_{i=1}^M \dlog z_i\in H^{M}(A-\cup_{i=1}^M \{z_i=0\},\Omega^M)$. Here the local generator of $\ml(a)\boxtimes \ml(-a)_0$ is the function $g g(0)^{-1}$ (the second coordinate is a constant), which is one plus a function that vanishes at zero.

To see that $g g(0)^{-1}\prod_{i=1}^M \dlog z_i$ produces a non-zero element in $H^{2M-1}(A-\{0\},\Bbb{C})$, we may apply
Grothendieck residues \cite[Page 651]{GH}, and hence the proof of  Theorem \ref{main} is complete. In fact we only need to show that $\prod \dlog z_i$ produces a non-zero element in $H^{2M-1}(A-\{0\},\Bbb{C})$, which is standard.

\section{Proof of Proposition 80}\label{keyF}
For the proof of the proposition assume, without loss of generality, that $q_1=1, q_2=2,\dots, q_{k}= k$. The left hand side of \eqref{formule1} is (with $(k-1)+b =M$, i.e., $k+b=M+1$),
\begin{multline}\label{formuleLeft}
\sum_{c=1}^{k}(-1)^{c+1} S_{1,\dots, \hat{c},\dots, k}\\
\ =\sum_{c=1}^{k}(-1)^{c+1} \prod_{q=1, q\neq c}^{k} \dlog(F^{(1)}_q-F^{(2)}_q)\bigg( \sum_{1\leq i_1<\dots<i_{b}\leq r}\prod_{s=1}^b a_{i_s}\dlog f^{(1)}_{i_s} \dlog f^{(2)}_{i_s}\bigg).
\end{multline}
The right hand side  of \eqref{formule1} is as follows
\begin{multline}\label{RHS}
(\eta^{(1)}-\eta^{(2)})\wedge S_{1\dots, k}\\
=(\eta^{(1)}-\eta^{(2)})\wedge \prod_{q=1}^{k}\dlog(F^{(1)}_q-F^{(2)}_q)\bigg(\sum_{ 1\leq j_1 <\dots < j_{b-1}\leq r} \prod_{s=1}^{b-1}{a_{j_s}}\dlog f^{(1)}_{j_s} \dlog f^{(2)}_{j_s}\bigg).
\end{multline}

Let us compare the coefficients of $a_{i_1}\dots a_{i_b}$ on both sides of \eqref{formule1}. Both coefficients are easily seen to be zero if $i_j=i_{j'}$ for some $j\neq j'$. We will assume that this is not the case.

Since $k+b=M+1$, we may write a linear dependence equation of form $\sum_{q=1}^{k} c_q F_{q} + \sum_{s=1}^b b_s f_{i_s}$ equals a constant with at least one coefficient not zero. We may assume that this constant is zero, and $c_1=1$ (if the $f_i$ are affinely dependent then both sides are zero). So after adjusting signs,
$$F_{1}=\sum_{q=2}^{k} c_q F_{q} + \sum_{s=1}^b b_s f_{i_s}.$$

We will now simply replace $F^{(1)}_{1}-F^{(2)}_{1}$ by
\begin{equation}\label{cat}
\sum_{q=2}^{k} c_q (F^{(1)}_{q}-F^{(2)}_q) + \sum_{s=1}^b b_s (f^{(1)}_{i_s}-f^{(2)}_{i_s})
 \end{equation}
 wherever it appears (which is all terms of RHS of \eqref{formule1}, and in all but one term on LHS of \eqref{formule1}).

The left hand side of \eqref{formule1} is
\begin{equation}\label{shopy}
\sum_{c=1}^{k}(-1)^{c+1} \prod_{q=1, q\neq c}^{k} \dlog(F^{(1)}_q-F^{(2)}_q)\prod_{s=1}^b \dlog f^{(1)}_{i_s} \dlog f^{(2)}_{i_s}
\end{equation}
which is a sum
\begin{equation}\label{shopy3}
\prod_{q=2}^{k} \dlog(F^{(1)}_q-F^{(2)}_q)\prod_{s=1}^b \dlog f^{(1)}_{i_s} \dlog f^{(2)}_{i_s}.
\end{equation}
plus sum over $c=2,\dots, k$ of
\begin{equation}\label{shopy4}
(-1)^{c+1} \dlog(F^{(1)}_1-F^{(2)}_1)\prod_{q=2,q\neq c}^{k} \dlog(F^{(1)}_q-F^{(2)}_q)\prod_{s=1}^b \dlog f^{(1)}_{i_s} \dlog f^{(2)}_{i_s}.
\end{equation}
The expression in \eqref{shopy4} equals (using \eqref{cat}) $$-\frac{\sum_{q=2}^{k}c_q(F^{(1)}_q-F^{(2)}_q)}{F^{(1)}_1-F^{(2)}_1}\cdot \prod_{q=2}^{k} \dlog(F^{(1)}_q-F^{(2)}_q)\prod_{s=1}^b \dlog f^{(2)}_{i_s} \dlog f^{(2)}_{i_s}$$

Therefore the coefficient of  $a_{i_1}\dots a_{i_b}$ on the LHS of \eqref{formule1} is
\begin{equation}\label{shopy5}
(1-\frac{\sum_{q=2}^{k}c_q(F^{(1)}_q-F^{(2)}_q)}{F^{(1)}_1-F^{(2)}_1})\cdot \prod_{q=2}^{k} \dlog(F^{(1)}_q-F^{(2)}_q)\prod_{s=1}^b \dlog f^{(1)}_{i_s} \dlog f^{(2)}_{i_s}.
\end{equation}

On the right hand side  of \eqref{formule1} the desired coefficient of  $a_{i_1}\dots a_{i_b}$ is sum over $\ell=1,\dots, b$ of
\begin{equation}\label{shopy2}
(\dlog f^{(1)}_{i_{\ell}} -\dlog f^{(2)}_{i_{\ell}})\wedge \bigg(\prod_{q=1}^{k}\dlog(F^{(1)}_q-F^{(2)}_q) \prod_{s=1, s\neq \ell}^{b} \dlog f^{(1)}_{j_s} \dlog f^{(2)}_{j_s}\bigg).
\end{equation}
We replace $\dlog(F^{(1)}_1-F^{(2)}_1)$ by \eqref{cat}
and calculate
\begin{equation}\label{wait}
(\dlog f-\dlog g)\wedge(df-dg) =(f-g)\dlog f \wedge \dlog g.
\end{equation}

The quantity \eqref{shopy2} is therefore equal to
 \begin{equation}\label{shopy9}
 \frac{\sum_{\ell=1}^b b_{\ell} (f^{(1)}_{i_{\ell}}-f^{(2)}_{i_{\ell}})}{F^{(1)}_1-F^{(2)}_1}\prod_{q=2}^{k}\dlog(F^{(1)}_q-F^{(2)}_q) \prod_{s=1}^{b} \dlog f^{(1)}_{j_s} \dlog f^{(2)}_{j_s}.
\end{equation}
It is now easy to see that \eqref{shopy5} and \eqref{shopy9} are equal and we are done.

\section{Applications to invariant theory}
We will use the notational set-up as in the introduction (Section \ref{IntroIn}).

\subsection{Conformal Blocks}\label{sanders} Let $\ell$ be a positive integer and consider an $n$-tuple of dominant integral weights $\vec{\lambda}=(\lambda_1,\dots,\lambda_n)$ such that each $(\lambda_i,\theta )\leq \ell)$ for each $1\leq i\leq n$. Let $\vec{z}=(z_1,\dots, z_n)$ be an $n$-tuple of distinct points of $\mathbb{A}^1 \subset \mathbb{P}^1$. Associated to this data there is the space of dual conformal blocks $\mathbb{V}_{\mathfrak{g},\vec{\lambda},\ell}(\mathbb{P}^1,z_1,\dots,z_n)$ which is a quotient of $\Bbb{A}(\vec{\lambda})$ (see the survey \cite{sorger}). The following is an explicit description of the quotient:

Define an operator
$$T_{\vec{z}}\in \op{End}(V(\vec{\lambda})),\  V(\vec{\lambda}) =V_{\lambda_1}\otimes \dots \otimes V_{\lambda_n}, \ T_{\vec{z}}=\sum_{i=1}^n z_i e_{\theta}^{(i)}$$
 with $e_{\theta}^{(i)}$ acting on the $i$-th tensor summand and $\theta$  the highest root. Let  $C_{\vec{z}}$ denote the image of  $T_{\vec{z}}^{\ell+1}$. The space of dual conformal blocks $\mathbb{V}_{\mathfrak{g},\vec{\lambda},\ell}(\mathbb{P}^1,z_1,\dots,z_n)$ is the cokernel of the natural map $C_{\vec{z}}$ to $\mathbb{A}(\vec{\lambda})$ \cite{Beau,FSV}. The dual vector spaces $\mathbb{V}^*_{\mathfrak{g},\vec{\lambda},\ell}(\mathbb{P}^1,z_1,\dots,z_n)$ are known as conformal blocks.

\subsection{The Schechtman-Varchenko map}\label{SVM}As in Section \ref{IntroIn}, let $V(\vec{\lambda})^*_0$ denote the zero $\mathfrak{h}$-weight space of $V(\vec\lambda)^*$.
Suppose $$\psi\in V(\vec{\lambda})^*_0=(V_{\lambda_1}\tensor V_{\lambda_2}\tensor\dots\tensor V_{\lambda_n})_0^*.$$ Let $v_i$ be a highest weight vector in $V_{\lambda_i}$, $i=1,\dots,n$.
Then, $${\Omega}^{SV}_{\vec{\lambda}}(\psi)=\Omega^{SV}(\psi)=\psi(v(\vec{t},\vec{z})) dt_1\dots dt_M\in A^M(U)$$ where for a fixed $\vec{z}$
$$v(\vec{t},\vec{z})=\sum_{\op{part}}\prod_{i=1}^n\langle\langle \prod_{b\in I_i} f_{\beta(b)}(t_b) v_i\rangle\rangle \in (V_{\lambda_1}\tensor V_{\lambda_2}\tensor\dots\tensor V_{\lambda_n})_0,$$
where $\beta:[1,\dots, M]\rightarrow R$ is a map to the positive roots as in Section \ref{IntroIn} and
$$\langle\langle f_{\gamma_1}(u_1)f_{\gamma_2}(u_2)\dots f_{\gamma_q}(u_q) v_i\rangle\rangle=\sum_{\op{perm}}\frac{1}{(u_1-u_2)(u_2-u_3)\dots(u_q-z_i)}
(f_{\gamma_1}f_{\gamma_2}\dots f_{\gamma_q}v_i).$$

Here, $\sum_{\op{part}}$ part stands for the summation over all partitions of $I = \{1,\dots, M\}$ into $n$
disjoint parts $I = I_1 \cup I_2 \cup · · · \cup I_n$ and $\sigma_{\op{perm}}$
perm the summation over all permutations of the
elements of $\{1,\dots, q\}$. The operators $f_{\gamma}$ are the standard $f$ operators corresponding to simple roots $\gamma$.
\subsection{Proof of Proposition 11}\footnote{We thank A. Varchenko for a useful discussion that led to the following proof.}
 Recall (see Section \ref{IntroIn}) that we need to show that
$$\mathbb{A}(\vec{\lambda})^*\hookrightarrow H^M(A^{\starr}(U), \eta\wedge).$$ The cohomology of the Aomoto complex  is independent of $\kappa$, and the map in Proposition \ref{injective} is linear in $\frac{1}{\kappa}$. Therefore an element in the kernel of the map for one value of $\kappa$, is in the kernel for any value of $\kappa$. We may therefore assume that $\kappa=\ell +g^*$ where $\ell$ is a sufficiently large integer, and $g^*$
is the dual Coxeter number. We take $\ell$ large so that $\Bbb{A}(\vec{\lambda})^*$ coincides with the space of conformal blocks $$\Bbb{V}^*_{\frg,\vec{\lambda},\ell}(\Bbb{P}^1,z_1,\dots,z_n)$$
for the data $\vec{\lambda}$ at level $\ell$ (with the marked curve equal to $\Bbb{P}^1$, and the marked points $z_1,\dots,z_n$).
Note that there is always a surjective map $\Bbb{A}(\vec{\lambda})\to \Bbb{V}_{\frg,\vec{\lambda},\ell}(\Bbb{P}^1,z_1,\dots,z_n)$.

It is a consequence of \cite{Ram,Bel}, and results of Deligne \cite{Del2} (see Section \ref{honda} below) that the map
\begin{equation}\label{trace}
\Bbb{V}^*_{\frg,\vec{\lambda},\ell}(\Bbb{P}^1,z_1,\dots,z_n)\to H^M(U,\ml(a))
\end{equation}
is injective (for large $\ell$). But \eqref{trace} factors through the map in Proposition \ref{injective}, and therefore Proposition \ref{injective} follows.

For any choice of weights, $H^M(A^{\starr}(U),\eta\wedge)=H^M(P,\underline{\ml}(a))$, and hence we have the following:
\begin{corollary}\label{satz}
The induced mapping $\Bbb{A}(\vec{\lambda})^*\to  H^M(P,\underline{\ml}(a))$ is injective for arbitrary $\kappa$.
\end{corollary}

\subsection{Some general notation}\label{honda}
We give a more detailed proof of the injectivity of the map \eqref{trace}, introducing notation that will be used elsewhere.

Let $C$ be a sufficiently divisible integer such that $C\kappa\eta$ is the logarithmic derivative $\frac{dQ}{Q}$ of a regular function $Q$ on $U$, here we may also assume that $\Sigma_M$ (see \eqref{symmetrie} preserves $Q$. Consider the (cyclic) unramified cover $\pi:\widehat{U}\to U$ defined by $\mathcal{R}^{C\kappa}=Q$. Note that

\begin{itemize}
\item $\Sigma_M\times \mu_{C\kappa}$ acts on $\widehat{U}$, where $\mu_{C\kappa}$ is the group of $C\kappa$ roots of unity in $\Bbb{C}^*$.
\item The multivalued function $\mathcal{R}$ on $U$ becomes single valued on $\widehat{U}$.

\item$\ml(a)$ is an isotypical component of $\pi_{*}\Bbb{C}$, for the character $\tau:\Sigma_M\times \mu_{C\kappa}\to  \Bbb{C}^*$ given by $\tau(\sigma,\mu)=\mu^{-1}\epsilon(\sigma)$ where $\epsilon$ is the sign character.
\end{itemize}

Let $\widehat{P}$ be a smooth
projective $\Sigma_M\times \mu_{C\kappa}$ equivariant compactification of $\widehat{U}$. It is shown in \cite{Ram,Bel} that the corresponding injective map
\begin{equation}\label{trace1}
\Bbb{V}^*_{\frg,\vec{\lambda},\ell}(\Bbb{P}^1,z_1,\dots,z_n)\to H^0(\widehat{U},\Omega^M)^{\tau}
\end{equation}
which sends $\psi\in \Bbb{V}^*_{\frg,\vec{\lambda},\ell}(\Bbb{P}^1,z_1,\dots,z_n)$ to $\mathcal{R}\Omega^{SV}(\psi)$
is induced by an injective map (i.e. the forms in the image of \eqref{trace1} extend to compactifications) $$\Bbb{V}^*_{\frg,\vec{\lambda},\ell}(\Bbb{P}^1,z_1,\dots,z_n)\to H^0(\widehat{P},\Omega^M)^{\tau}=H^{M,0}(\widehat{P},\Bbb{C})^{\tau}.$$

By results of Deligne (see \cite{Ram}), the map $H^{M,0}(\widehat{P})\to H^M(\widehat{U},\Bbb{C})$ is injective. Therefore the map
\begin{equation}\label{inje}
\Bbb{V}^*_{\frg,\vec{\lambda},\ell}(\Bbb{P}^1,z_1,\dots,z_n)\to H^M(\widehat{U},\Bbb{C})^{\tau}=H^M(U,\ml(a))
\end{equation}
is injective. Therefore \eqref{trace} is injective.

\subsection{Proof of a generalization of Theorem 20}
We consider the maps \eqref{SVmap} for two sets of data:
\begin{enumerate}
\item The representations $(\lambda_1,\dots,\lambda_n)$, and $\kappa$.
\item The representations $(\lambda^*_1,\dots,\lambda^*_n)$, and $-\kappa$.
\end{enumerate}
As in the introduction (see Section \ref{IntroIn}), these two give rise to the same hyperplane arrangement, but with weights that are negatives of each other. We therefore find a compactification $P\supseteq U$, and two objects (see Section \ref{mostgeneralform}) in the derived category $D^b_c(P)$: $\underline{\ml}(a)$ and $\underline{\ml}(-a)$.

The map \eqref{beforeclass}, and Proposition \ref{injective} gives rise to the following two injective maps (it should really be $-\eta$ in the second equation but the sign does not affect the quotient):
\begin{equation}\label{SVmap1}
\Bbb{A}(\vec{\lambda})^{*}\hookrightarrow \bigl(\frac{A^M(U)}{\eta\wedge A^{M-1}(U)}\bigr)^{\chi}.
\end{equation}
and
\begin{equation}\label{SVmap2}
\Bbb{A}(\vec{\lambda}^*)^*\hookrightarrow\bigl(\frac{A^M(U)}{\eta\wedge A^{M-1}(U)}\bigr)^{\chi}.
\end{equation}
Recall the form $S \in A^{2M}(U\times U)$ as in the statement of Theorem \ref{maino}. We will use $S$ to form a diagram, which will be shown to commute:
\begin{equation}\label{RedSquare}
\xymatrix{
(\bigl(\frac{A^M(U)}{\eta\wedge A^{M-1}(U)}\bigr)^{\chi})^*\ar[d]^{\tilde{S}} \ar@{->>}[r] & \Bbb{A}(\vec{\lambda}^*)\ar[d]^{\sim}\\
\bigl(\frac{A^M(U)}{\eta\wedge A^{M-1}(U)}\bigr)^{\chi}   &          \Bbb{A}(\vec{\lambda})^*.\ar@{_{(}->}[l]
}
\end{equation}

For a representation $V$ of the group $\Sigma_M$, the natural map
\begin{equation}\label{nature}
(V^*)^{\chi}\to (V^{\chi})^*
\end{equation}
is an isomorphism. The map $\tilde{S}$ in \eqref{RedSquare} arises as follows. From $S$, one obtains a $\Sigma_M$-equivariant map $A^M(U)^*\to A^M(U)$, taking $\chi$- isotypical components, we get $(A^M(U)^*)^{\chi}\to A^M(U)^{\chi}$. Composing with the inverse of the natural map \eqref{nature}, we get the map $\tilde{S}$ in \eqref{RedSquare}.
The  vertical map on the right of \eqref{RedSquare} is the inverse of the natural isomorphism \eqref{canaan} from invariants $\mathbb{A}(\vec{\lambda})^*$ to coinvariants $\mathbb{A}(\vec{\lambda}^*)$.

The following result is a direct  consequence of \cite[Theorem 6.6]{SV}, as we will explain in the next section.
\begin{proposition}\label{WSquare}
The diagram \eqref{RedSquare} commutes up to  a non-zero multiplicative constant.
\end{proposition}
Putting together the $\chi$ isotypical component of \eqref{zhGCC} together with \eqref{RedSquare}, we get a diagram which commutes up to  a non-zero factors (commutes projectively):
\begin{equation}\label{FitBit}
\xymatrix{
H^M(P,{D}(\underline{\ml}(-a))[-2M])^{\chi}\ar[dr]^{\alpha} \ar[r]^{\sim} & (H^{M}(P, \underline{\ml}(-a))^*)^{\chi}\ar[d]^{\Sigma}\ar[r]^{\sim} & (\bigl(\frac{A^M(U)}{\eta\wedge A^{M-1}(U)}\bigr)^{\chi})^*\ar[d]^{S} \ar@{->>}[r] & \Bbb{A}(\vec{\lambda}^*)\ar[d]^{\sim}\\
& H^M(P,\underline{\ml}(a))^{\chi} & (\bigl(\frac{A^M(U)}{\eta\wedge A^{M-1}(U)}\bigr)^{\chi}) \ar[l]^{\sim} &          \Bbb{A}(\vec{\lambda})^*.\ar@{_{(}->}[l]
}
\end{equation}
This leads to the following description, of the image of the map in Corollary \ref{satz}:
\begin{theorem}\label{egregium}
The image of the injective mapping $\Bbb{A}(\vec{\lambda})^*\to  H^M(P,\underline{\ml}(a))^{\chi}$ coincides with the image of the
topological map $H^M(P,D(\underline{\ml}(-a))[-2M])^{\chi}\to  H^M(P,\underline{\ml}(a))^{\chi}$ (see the map \eqref{zhGC}). Therefore,
\begin{equation}\label{MP}
\Bbb{A}(\vec{\lambda})^*=\operatorname{Image}: H^M(P,D(\underline{\ml}(-a))[-2M])^{\chi}\to  H^M(P,\underline{\ml}(a))^{\chi}.
\end{equation}
\end{theorem}
Theorem \ref{egregium}  specializes to Theorem \ref{secondmain} for large $|\kappa|$.
\section{Proof of Proposition 106}
There is another way of obtaining the Schechtman-Varchenko maps \eqref{SVmap} with the role of $e$'s replaced by $f$'s and the highest weight vectors replaced by lowest weight vectors:
Let $$\psi\in V(\vec{\lambda})^*_0=(V_{\lambda_1}\tensor V_{\lambda_2}\tensor\dots\tensor V_{\lambda_n})_0^*.$$
 The lowest weight vector in $V_{\lambda_i}$ is $w_0(v_i)$ where $v_i$ is the highest weight vector in $V_{\lambda_i}$, here $w_0$ is a lifting of the longest element of the Weyl group to $G$. Now, the weight of $w_0(v_i)$ is $w_0(\lambda_i)$ and
$$\sum_{i=1}^n w_0(\lambda_i)= \sum_{b=1}^n w_0(\beta(b)).$$
We can form a map
$$\Omega^{SV-}: V(\vec{\lambda})_0^*\to A^M(U)$$ again
with $e_i$ replaced by $f_{i'}$ where $\alpha_{i'}=w_0(\alpha_i)$.
The following is easy to check:
\begin{lemma}\label{debate}
$$\Omega^{SV-}(\psi)=\Omega^{SV}(w_0(\psi))$$
and hence if $\psi$ is Weyl group invariant then $\Omega^{SV-}(\psi)=\Omega^{SV}(\psi)$.
\end{lemma}
Let $$\op{Sh}:V(\vec{\lambda}^*)_0^*\to V(\vec{\lambda})^*_0$$
be the Shapovolov isomorphism in representation theory induced from an isomorphism without the outer duals and weight spaces. Let $\psi\in V(\vec{\lambda^*})_0^*$:
then it follows that $$\Omega^{SV-}_{\vec{\lambda}^*}(\psi)=\Omega^{SV}_{\vec{\lambda}}(\op{Sh}(\psi)).$$

\subsection{}
We will use $\Omega^{SV-}$ for the Schechtman-Varchenko map for $(\lambda_1^*,\dots,\lambda_n^*)$. This does not change \eqref{WSquare} since  we have restricted the maps \eqref{SVmap} always to invariants in forming that diagram, and Lemma \ref{debate}. We will also identify $V(\vec{\lambda}^*)_0^*$ and  $V(\vec{\lambda})^*_0$, as well as  $\Bbb{A}(\vec{\lambda}^*)$ and $\Bbb{A}(\vec{\lambda})$ by the Shapovolov form.

Consider the following diagram. Here the objects on the top row are to be considered as objects for the dual weights, using the above identification. Therefore the object on the top right is also $\Bbb{A}(\vec{\lambda}^*)$. The vertical map
$\Bbb{A}(\vec{\lambda})\to \Bbb{A}(\vec{\lambda})^*$ is the inverse of the evident isomorphism $\Bbb{A}(\vec{\lambda})^*\to \Bbb{A}(\vec{\lambda})$.

\begin{equation}\label{OrangeSquare1}
\xymatrix{
(A^M(U)^{\chi})^*\ar[d]^{\tilde{S}}\ar[r]^{{\Omega^{SV}_{\vec{\lambda}}}^*}&  V(\vec{\lambda})_0\ar[d]^{\op{Sh}}\ar[r] & \Bbb{A}(\vec{\lambda})\ar[d]\\
{A^M(U)}^{\chi}   &         V(\vec{\lambda})_0^*\ar[l]^{\op{SV}_{\lambda}} & \Bbb{A}(\vec{\lambda})^*.\ar[l]
}
\end{equation}
\begin{enumerate}
\item The square on the right almost commutes, i.e., the two maps $V(\vec{\lambda})_0\to  V(\vec{\lambda})^*_0$ are not the same. They differ by a map to $(n_{+}V(\vec{\lambda})^*)_0$.
\item The square on the left commutes: This is  \cite[Theorem 6.6]{SV}, noting that the flag complex of degree $M$ in \cite{SV} is dual to the space of logarithmic forms $A^M(U)$ (see Theorem 2.4 in loc. cit, and below).
\end{enumerate}
The almost commutativity of the big square resulting from the above picture  above allows us to conclude that \eqref{WSquare} commutes: The  map \eqref{SVmap} is $V(\vec{\lambda})_0^*\to A^M(U)$ has the following property (see \cite[Lemma 6.8.10]{SV}) $$(n_{+}V^*(\vec{\lambda}))_0\mapsto\eta\wedge A^{M-1}(U).$$

\subsection{Flag forms}We now recall the notion of flag forms \cite[Section 3]{SV} $\mathcal{F}^p$ for $0\leq p\leq M$. Let as before $\mathcal{C}$ be an arrangement given by linear polynomials $f_1,\dots, f_r$ and let $H_i=Z(f_i)$ and  $\mathcal{C}^i$ denote the set of edges of codimension $i$ in $\mathbb{A}^M$.  For $0\leq p\leq M$, let $\widetilde{\mathcal{F}}^p$ be the set of flags $(L^0\supset L^1\dots \supset L^p)$, where $L^i \in \mathcal{C}^i$. Let $\mathcal{F}^p$ be the quotient of the free abelian group generated by the set of flags $\widetilde{\mathcal{F}}^p$ by relations of the following form:
$$\sum_{F\supset \widehat{F}}F=0,$$ where $\widehat{F}=(L^0\supset \dots L^{i-1} \supset L^{i+1}\supset \dots \supset L^p)$ and $F=(\tilde{L}^0\supset \dots \supset \tilde{L}^p)$ such that $\tilde{L}^j=L^j$ for all $j\neq i$. 

Now for each $p$, following \cite{SV}, we define a map $\varphi^p: {A}^p(U)\rightarrow {\mathcal{F}^p}^*$. For a $p$-tuple of hyperplanes $(H_1,\dots, H_p)$ in general position,  we define $F(H_1,\dots, H_p)=(H_1\supset H_{12}\supset H_{123}\supset \dots \supset H_{12\dots p}) \in \mathcal{F}^p$, where $H_{1\dots i}=H_1\cap\dots\cap H_i$. These flag forms are dual to the Aomoto log forms \cite[Theorem 2.4]{SV} by the map $\varphi^p$ is given by the following formula

\begin{equation}\label{svdual}
\varphi^p(H_1,\dots, H_p)=\sum_{\sigma \in S_p}(-1)^{sgn(\sigma)}\delta_{F(H_\sigma(1),\dots H_{\sigma(p)})},
\end{equation}
where $\delta_F$ for a flag $F$ is the delta functional for a flag $F$.

\subsubsection{Quasi-classical contravariant form}Consider the arrangement $\mathcal{C}$ and $a$ be a set of weights. Let $\bar{H}=(H_1,\dots, H_p) \in A^p(U)$, we say $H$ is adjacent to a flag $F$ if we can find a permutation $\sigma \in S_p$ such that $F=F(H_{\sigma(1)},\dots, H_{\sigma(p)})$. For any two flags $F$ and $F'$, Schechtman-Varchenko (\cite[Equation (3.3.1)]{SV}) defines a bilinear form on $\mathcal{F}^p$.
\begin{equation}\label{SVform}
S^p(F,F')=\frac{1}{p!}\sum_{\bar{H}}(-1)^{\sigma_{\bar{H}}(F)\sigma_{\bar{H}}(F')}a(H_1)\dots a(H_p),
\end{equation} where the sum is taken over all $\bar{H}=(H_1,\dots, H_p)$ adjacent to the flags $F$ and $F'$.
We will use the following fact that is an easy check using Equation \eqref{SVform}

\begin{proposition}
Under this identification, the quasi-classical contravariant form on $\mathcal{F}^M$ equals the
Shapovolov form $S^{(M)}$ given by \eqref{diagonalgeneral}.
\end{proposition}

\begin{proof}

This proceeds as follows: Pick two flags $F$ and $F'$: The two values on the pair $(F,F')$ are sums over $(H_1,\dots, H_M)$. Start with the values assigned by our $S$: a term $$\Omega=(H^{(1)}_1H_1^{(2)}\dots H_M^{(1)}H_M^{(2)})$$ acts on $F\tensor F'$ with non-zero coefficient only if $F$ is of the form
$$F=(H_{\sigma(1)}\supset H_{\sigma(1)}\cap H_{\sigma(2)}\supset \dots ),$$ similarly $F'$ for a permutation $\sigma$. In the formula for  $\Omega$, $H^{(1)}_i= \dlog f_i^{(1)}$ where $f_i$ is the defining equation for $H_i$, similarly for  $H^{(2)}_i$. This coefficient is so  because a log form acts on flag forms via the formula \cite[Equation 2.3.2]{SV}. The corresponding coefficient is the product of the weights $a(H_i)$ times $(-1)$ raised to the sum of signs of $\sigma$ and $\sigma'$.

For the form defined in \cite[Equation (3.3.1)]{SV}, the contributions are again over the same $\sigma $ and $\sigma'$, and the signs are the same.
\end{proof}

\subsection{Aomoto  Representatives}\label{labelo}
Consider the image $\Omega_{\psi}$ of the map \eqref{SVmap} on an element $\psi\in \Bbb{A}(\vec{\lambda})^*=(V(\vec{\lambda})^*)^\frg$. We claim that $\Omega_{\psi}$ does not have poles on any $E_{\alpha}$ with $a_{\alpha}=0$. This is because we may raise the level (and assume it to be integral) so that $\psi$ lives in the space of conformal blocks. Then (see Section \ref{honda} for the notation), the form $\mathcal{R}\Omega_{\psi}$ extends to compactifications of ramified covers of $P$ by the main results of \cite{Ram,Bel}, and is therefore square integrable on $P$. Thus on any $E_{\alpha}$ on which $\mathcal{R}$ has order zero, the form $\Omega_{\psi}$ is regular. The claim follows immediately because $a_{\alpha}$ is the order of $\mathcal{R}$ along $E_{\alpha}$.

The image of the mapping $\Bbb{A}(\vec{\lambda})^*\to  H^M(P,\underline{\ml}(a))^{\chi}$, which coincides with the image of  $H^M(P,D(\underline{\ml}(-a))[-2M])^{\chi}\to  H^M(P,\underline{\ml}(a))^{\chi}$ (by Theorem \ref{egregium}), is therefore contained in the span of $[\Omega]$ where $\Omega$ runs through elements of $A^M(U)^{\chi}$ which do not have poles on divisors $E_{\alpha}$ with $a_{\alpha}=0$. Such elements  $[\Omega]$ are always in the image of the topological mapping $H^M(P,D(\underline{\ml}(-a))[-2M])^{\chi}\to  H^M(P,\underline{\ml}(a))^{\chi}$ (see Lemma \ref{linearqQ}). Therefore,
\begin{proposition}\label{reverso}
The image of the Schechtman-Varchenko  mapping \eqref{SVmap}:
$$\Bbb{A}(\vec{\lambda})^*\to  H^M(P,\underline{\ml}(a))^{\chi}.$$
 equals the  set of $[\Omega]$ where $\Omega$ runs through elements of $A^M(U)^{\chi}$ which do not have poles on divisors $E_{\alpha}$ with $a_{\alpha}=0$.
\end{proposition}
The classes $[\Omega]$ in this proposition can be zero even if $\Omega$ is not zero.
\subsection{Mixed Hodge structures}\label{mathew}
Assume that $\kappa\neq 0$ is an integer. Recall \eqref{MP}:
$$\Bbb{A}(\vec{\lambda})^*=\operatorname{Image}: H^M(P,D(\underline{\ml}(-a))[-2M])^{\chi}\to  H^M(P,\underline{\ml}(a))^{\chi}.$$
Note that as in Section \ref{Lgen},  we may write $$ H^M(P,\underline{\ml}(a))^{\chi}=H^M(V,j_{!}\mathcal{L}(a))^{\chi}, \ H^M(P,D(\underline{\ml}(-a))[-2M])^{\chi}= H^M(V',q_{!}\ml(a)).$$

Using notation from Section \ref{honda}, it is easy to see that $H^M(V,j_{!}\mathcal{L}(a))^{\chi}$ is an isotypical component $H^M(\widehat{V},D;\Bbb{C})^{\tau}$ where $D\subset \widehat{V}$ is the inverse image $\pi^{-1}(V-U)$ (a closed subset), and hence acquires a mixed Hodge structure over the cyclotomic field $\Bbb{Q}(\mu_{C\kappa})$; here $\tau:\Sigma_M\times \mu_{C\kappa}\to \Bbb{C}^*$ is a character.

Similarly  $H^M(V',q_{!}\mathcal{L}(a))^{\chi}$ is isomorphic to $H^M(\widehat{V'},D';\Bbb{C})^{\tau}$ where $D'\subset \widehat{V'}=\pi^{-1}(V')$ is the inverse image $\pi^{-1}(V'-U)$, and hence acquires a mixed Hodge structure (MHS) over a cyclotomic field. Therefore,
\begin{itemize}
\item By Theorem \ref{egregium}, $\Bbb{A}(\vec{\lambda})^*$ the image of a morphism of mixed Hodge structures,  acquires a mixed Hodge structure defined  over a cyclotomic field.
 \end{itemize}
 Section \ref{PS} gives an example where this mixed Hodge structure is not pure.

There is an exact sequence of topological cohomology groups with $\Bbb{C}$ coefficients,
 $$H^{M-1}(D)\to H^M(\widehat{V},D)\to H^M(\widehat{V})\to H^M(D).$$
In the following, $F^{\starr}$ is the Hodge filtration and $M=\dim V$. We record some facts about the Hodge filtration
 \begin{enumerate}
 \item[(i)] $F^{M+i}(H^{M-1}(D))=0$,  $F^{M+i}(H^{M}(D))=0$, $i=0,1$, and hence exactness of the functors $F^k$ implies
 \item[(ii)] $F^{M+i}(H^{M}(\widehat{V},D))=F^{M+i}( H^M(\widehat{V})),i=0,1$.
 \item[(iii)] $F^M( H^M(\widehat{V}))\cap W_M( H^M(\widehat{V}))=H^{M,0}(\widehat{P})$ and  $F^{M+1}( H^{M}(\widehat{V}))=0$.

 (Note that $W_M( H^M(\widehat{V}))$ is the image
 of $H^{M}(\widehat{P})$ in $H^{M}(\widehat{V})$, and $H^{M,0}(\widehat{P})$ injects into $H^M(\widehat{V})$ as global $M$-forms on $\widehat{P}$ are also log forms on $\widehat{V}$.)
 \item[(iv)] $F^{M+1}H^M(\widehat{V},D)=0$.
 \end{enumerate}
 For (i), by \cite[page 45 (e)]{Del3}, $h^{p,q}$ for $H^n(D)$ vanishes if $p>\dim D =M-1$.

The above applies also to $(\widehat{V'},D')$ and hence both $H^M(\widehat{V'},D';\Bbb{C})^{\tau}$ and $H^M(\widehat{V},D;\Bbb{C})^{\tau}$ have $F^M\cap W_M$ equal to $H^{M,0}(\widehat{P})^{\tau}$. Therefore,
\begin{proposition}
The above MHS on  $\Bbb{A}(\vec{\lambda})^*$ has $F^M\cap W_M$  isomorphic to  $H^{M,0}(\widehat{P})^{\tau}$.
\end{proposition}
\subsection{}
Now consider the case $\kappa=\ell + g^*$ with $\ell$ a positive integer. The space of invariants $\Bbb{A}(\vec{\lambda})^*$ has a subspace given by the space of conformal blocks
 $\Bbb{V}^*_{\frg,\vec{\lambda},\ell}(\Bbb{P}^1,z_1,\dots,z_n)$. 	
 The subspace of conformal blocks $\Bbb{V}^*_{\frg,\vec{\lambda},\ell}(\Bbb{P}^1,z_1,\dots,z_n)$ injects into  $H^{M,0}(\widehat{P})^{\tau}$
 by \cite{Ram,Bel}. The image of conformal blocks is all  of $H^{M,0}(\widehat{P})^{\tau}$ for classical groups and $G_2$ \cite{BM}. Therefore
\begin{proposition}
 For classical groups and $G_2$,
\begin{equation}\label{nice}
W_M(\Bbb{A}(\vec{\lambda})^*)\cap F^M(\Bbb{A}(\vec{\lambda})^*)= \Bbb{V}^*_{\frg,\vec{\lambda},\ell}(\Bbb{P}^1,z_1,\dots,z_n)
\end{equation}
\end{proposition}

For $\ell$ sufficiently large, $\Bbb{A}(\vec{\lambda})^*$ coincides with $\Bbb{V}^*_{\frg,\vec{\lambda},\ell}(\Bbb{P}^1,z_1,\dots,z_n)$  and therefore $\Bbb{A}(\vec{\lambda})^*$ carries a pure Hodge structure of type $(M,0)$ over a cyclotomic field.

 \begin{remark}
If $\ell>-1+ \frac{1}{2}(\lambda_i,\theta)$, then $\Bbb{A}(\vec{\lambda})^*=\Bbb{V}^*_{\frg,\vec{\lambda},\ell}(\Bbb{P}^1,z_1,\dots,z_n)$ (see \cite{BGM} where this is derived from \cite{Beau,FSV}).
 \end{remark}
\subsection{Hodge filtration}
By the previous section, for $\frg$ classical or the Lie algebra of the exceptional group $G_2$,
$$W_M(\Bbb{A}(\vec{\lambda})^*)\cap F^M(\Bbb{A}(\vec{\lambda})^*)=\Bbb{V}^*_{\frg,\vec{\lambda},\ell}(\Bbb{P}^1,z_1,\dots,z_n).$$

The following seems to be plausible:
\begin{question}\label{plausible}
\item[(1)] Is $W_M(\Bbb{A}(\vec{\lambda})^*)=\Bbb{V}^*_{\frg,\vec{\lambda},\ell}(\Bbb{P}^1,z_1,\dots,z_n)$, and hence $W_{M-1}(\Bbb{A}(\vec{\lambda})^*)=0$.
\item[(2)] Assuming (1), consider the induced MHS on $T=\Bbb{A}(\vec{\lambda})^*/W_M(\Bbb{A}(\vec{\lambda})^*)$. For $p\geq 0$, is
 $F^{M-p}(T)=\Bbb{V}^*_{\frg,\vec{\lambda},\ell+p+1}(\Bbb{P}^1,z_1,\dots,z_n)/W_M(\Bbb{A}(\vec{\lambda})^*)$?
\end{question}g
The only justification we have for (2) is that the KZ connection behaves consistently as we see below.

By the description of conformal blocks given in Section \ref{sanders},
$$\Bbb{V}_{\ell}\mid_{\vec{z}}=
\Bbb{V}_{\frg,\vec{\lambda},\ell}(\Bbb{P}^1,z_1,\dots,z_n)= \frac{V(\vec{\lambda})}{\frg V(\vec{\lambda}) +\im T_{\vec{z}}^{\ell+1}}.$$

\subsubsection{Connections}For each $\ell$, the KZ connection is a connection $\nabla^{(\ell)}$ on $(V_{\lambda_1}\tensor\dots\tensor V_{\lambda_n})$ of the form (see \eqref{wakeup})
$\frac{\partial}{\partial z_i}\mapsto \nabla^{(\ell)}_i=\frac{\partial}{\partial z_i} +\frac{A_i}{\ell+g^*}$ where
$A_i=A_i(\vec{z})$ are operators on $V(\vec{\lambda})$ which are independent of $\ell$.

It is known that $\nabla^{(\ell)}$ descends to a connection operator on $\Bbb{A}(\vec{\lambda})$ and on its  quotient
$\Bbb{V}_{\ell}$. Therefore $\frac{\partial}{\partial z_i} +\frac{A_i}{\ell+g^*}$ preserves the subbundle
$\frg (V(\vec{\lambda})) +\im T_{\vec{z}}^{\ell+1}$ of $V(\vec{\lambda})$.

This means that
\begin{lemma}
$[\frac{\partial}{\partial z_i} +\frac{A_i}{\ell+g^*}](T_{\vec{z}}^{\ell+1}w)$ is in $\frg V(\vec{\lambda}) +\im T_{\vec{z}}^{\ell+1}$.
\end{lemma}

But $\frac{\partial}{\partial z_i} T_{\vec{z}}^{\ell+1} w$ is a multiple of  $z_i T_{\vec{z}}^{\ell} e^{(i)}_{\theta}w$. Therefore the result of the  action of ${A_i}$ on the vector $T_{\vec{z}}^{\ell+1}w$
is in $\frg V(\vec{\lambda}) +\im T_{\vec{z}}^{\ell+1}$ up to  a multiple of $z_i T_{\vec{z}}^{\ell} e^{(i)}_{\theta}w$.
Therefore
\begin{proposition}
\begin{enumerate}
\item $\nabla^{(k)}$ carries the subbundle $\frg V(\vec{\lambda}) +\im T_{\vec{z}}^{\ell+1}$ of $ V(\vec{\lambda})$ to  $\frg V(\vec{\lambda}) +\im T_{\vec{z}}^{\ell}$ for any $k>0$.
\item Dualizing, the connection $\nabla^{(k)}$ on $\Bbb{A}(\lambda)^*$ carries subbundle $\Bbb{V}^*_{\ell}$ to the ``larger" $\Bbb{V}^*_{\ell+1}$ for any $\ell$. For $k=\ell$, $\nabla^{(k)}$ carries $\Bbb{V}^*_{k}$ to itself.
\item The filtration of $\Bbb{A}(\vec{\lambda})^*$ by the filtration of conformal blocks behaves like the Hodge filtration with respect to connections (Griffiths transversality).
\end{enumerate}
\end{proposition}

\section{Non semi-simple monodromy for $\mathfrak{sl}(2)$}\label{PS}

Consider the representation $\mathbb{C}^2$ of $\mathfrak{sl}(2)$ with highest weight $\omega_1$. Let $\vec{\lambda}=(\omega_1,\omega_1,\omega_1,\omega_1)$, and consider the space of coinvariants $\mathbb{A}(\vec{\lambda})$. Let $\vec{z}=(z_1,\dots,z_4)\in \Bbb{C}^4$ be a tuple of distinct points. As $\vec{z}$ varies in the configuration space of points in $\mathbb{A}^1$, we know that the trivial vector bundle $\mathcal{A}(\vec{\lambda})$ with fibers $\mathbb{A}(\vec{\lambda})$ is endowed with a flat connection known as  the Knizhnik-Zamolodchikov (KZ) connection. In this section, we show that the monodromy of the KZ connection on $\mathbb{A}(\vec{\lambda})$ with $\kappa=3$ is not semi-simple.

Consider the standard basis $v_1$ and $v_2$ of $\mathbb{C}^2$. We know that $\dim \mathbb{A}(\vec{\lambda})$ is two. Let $v=v_1\otimes v_1\otimes v_2 \otimes v_2$ and $w=v_1\otimes v_2\otimes v_1 \otimes v_2$. It is easy to check that the classes $[v]$ and $[w]$ in $\mathbb{A}(\vec{\lambda})$ form a basis. Let $e$, $f$ and $h$ be standard elements of $\mathfrak{sl}(2)$ and let $\Omega=\frac{1}{2}h^2 + e.f +f.e$ denote the Casimir operator considered as an element of the universal enveloping algebra of $\mathfrak{sl}(2)$.

Associated to the data $\vec{\lambda}$, the equations of the flat sections of KZ connections as operators on $\mathbb{A}(\vec{\lambda})$ are given by the following connection equations: Here $u\in \mathcal{A}(\vec{\lambda})=\Bbb{A}(\vec{\lambda})\tensor \mathcal{O}$,

\begin{equation}\label{wakeup}
\bigl(\frac{\partial}{\partial z_j}+\frac{1}{\kappa}\sum_{k=1,k\neq j}^4\frac{\Omega_{jk}}{z_j-z_k}\bigr)u =0,
\end{equation}
We now compute the matrices $\Omega_{jk}$
where $\kappa$ is a complex number, $\Omega_{jk}$ is the Casimir operator acting on $j$-th and $k$-th component with respect to the chosen basis $\{[v],[w]\}$,. They are given as follows:
\[ \Omega_{1,2}=\left( \begin{array}{cc}
1/2 & -1 \\
0 & -3/2
\end{array} \right), \  \Omega_{1,3}=\left( \begin{array}{cc}
-3/2 & 0 \\
-1 & 1/2
\end{array} \right), \  \Omega_{1,4}=\left( \begin{array}{cc}
-1/2 & 1 \\
1 & -1/2
\end{array} \right)\]

\[\Omega_{2,3}= \left( \begin{array}{cc}
-1/2 & 1  \\
1 & -1/2  \end{array} \right), \ \Omega_{2,4}=\left( \begin{array}{cc}
-3/2 & 0  \\
-1 & 1/2
\end{array} \right), \ \Omega_{3,4}=
\left( \begin{array}{cc}
1/2 & -1  \\
0 & -3/2  \end{array}\right). \]

There is a non-constant  rank one sub-bundle $\mathcal{K}(\vec{\lambda})$ of $\mathcal{A}(\vec{\lambda})$ given fiberwise by the kernel of the map $\mathbb{A}(\vec{\lambda})$ to the dual conformal blocks $\mathbb{V}_{\frg, \vec{\lambda}, 1}(\mathbb{P}^1,{z_1,\dots,z_n})$. Since the KZ connection on the bundle $\mathcal{A}(\vec{\lambda})$ descends to the bundle of conformal blocks, it follows that the bundle $\mathcal{K}$ is preserved by the KZ connection. Hence the monodromy representation of KZ connection of $\mathbb{A}(\vec{\lambda})$ is reducible.  We can describe the dual conformal blocks $\mathbb{V}_{\frg, \vec{\lambda}, 1}(\mathbb{P}^1,{z_1,\dots,z_n})$ as the class of $[v]$ or $[w]$ with the following relation:
\begin{equation}
[v]=-\frac{(z_1-z_2)(z_3-z_4)}{(z_1-z_3)(z_2-z_4)}[w].
\end{equation}
Further the fiber of $\mathcal{K}$ at $\vec{z}$ is given by  $\Phi(\vec{z})=(z_1-z_3)(z_2-z_4)v+(z_1-z_2)(z_3-z_4)w$ as an element of $\mathbb{A}(\vec{\lambda})$.  The following lemma can be checked by a direct calculation
\begin{lemma}\label{aflat}
With the above notations, we get
\begin{enumerate}
\item The multi-valued section of $\widetilde\Phi(\vec{z})=\prod_{1\leq i<j\leq 4}(z_i-z_j)^{-\frac{1}{6}}\Phi(\vec{z})$ is flat for the KZ connection on $\mathcal{A}(\vec{\lambda})$.
\item The multi-valued section $\widetilde{[v]}:=f[v]$ gives a flat section of the KZ connection in the dual conformal blocks bundle, where
 $$f(\vec{z})=(z_1-z_2)^{-1/2}(z_1-z_3)^{1/2}(z_1-z_4)^{1/2}(z_2-z_3)^{1/2}(z_2-z_4)^{1/2}(z_3-z_4)^{-1/2}.$$
\item The sections $\widetilde{v}=f(\vec{z}).v$ and $\widetilde{\Phi}(\vec{z})$ are related by the following:
$$\nabla_{\frac{\partial}{\partial z_j}}\widetilde{v}=-\frac{1}{3}\prod_{k=1,k\neq j}^4(z_j-z_k)^{-1}\prod_{1\leq i<k\leq 4}(z_i-z_k)^{\frac{1}{6}}.f(\vec{z}).\widetilde{\Phi}(\vec{z}),$$ where $f(\vec{z})$ is as above.
\end{enumerate}
\end{lemma}
The following proposition is the main result of this section
\begin{proposition}\label{nonsplit}
The monodromy representation of the KZ connection on $\mathbb{A}(\vec{\lambda})$ with $\kappa=3$ is not semisimple.
\end{proposition}
\begin{proof}We know that the monodromy representation is reducible and two dimensional. Assume that the representation is semi-simple, and hence abelian. We choose the point $(-1/2,0, 1/2,1)$ as our base point  and consider a Pochhammer loop $\gamma$, where the point $z_1$ moves and the other coordinates remain fixed.

We know that the Pochhammer loop is an element of the commutator of the fundamental group. Now if the monodromy is abelian, then the image of $\gamma$ under the monodromy representation must be the $2\times 2$-identity matrix. We have already shown that the monodromy is upper triangular. Hence the Pochhammer contour $\gamma$ maps to an unipotent matrix of the form
$$\left( \begin{array}{cc}
1 & 0  \\
a_{2,1} & 1  \end{array} \right).
$$
If we can show that the constant $a_{2,1}$ is non-zero, we will be done. We need to find a flat section $\Psi(\vec{z})$ of the form $\widetilde{v}+g(\vec{z})\widetilde{\Phi}(\vec{z})$. If such a section exists, then by Lemma \ref{aflat}, we will get that
$$\nabla_{\frac{\partial}{\partial z_1}} \Psi(\vec{z})=\xi(\vec{z})\widetilde{\Phi}(\vec{z})+ \frac{\partial}{\partial z_1}g(\vec{z}) \widetilde{\Phi}(\vec{z})=0,\ \mbox{i.e} \ \frac{\partial}{\partial z_1}g(\vec{z})=- \xi(\vec{z})$$
where $${\xi}(\vec{z})=-\frac{1}{3}(z_1-z_2)^{-4/3}(z_1-z_3)^{-1/3}(z_1-z_4)^{-1/3}(z_2-z_3)^{2/3}(z_2-z_4)^{2/3}(z_3-z_4)^{-1/3}.$$ The entry $a_{21}$ is exactly the  monodromy of  $g(\vec{z})$ about the Pochhammer contour $\gamma$. Since we are interested in the monodromy as $z_1$ changes, the terms $(z_2-z_3)^{2/3}(z_2-z_4)^{2/3}(z_3-z_4)^{-1/3}$ in $\xi(\vec{z})$ have no contribution to monodromy. Thus, we can just focus on the term $(z_1-z_2)^{-4/3}(z_1-z_3)^{-1/3}(z_1-z_4)^{-1/3}$. Let ${}_2F_{1}(a,b,c|u)$ denote the Gauss hypergeometric function.

Now for any $b,b-c \notin \mathbb{Z}$, we have
$${}_2F_{1}(a,b,c|u)=\frac{\Gamma(c)}{\Gamma(b)\Gamma (c-b)}\frac{1}{(1-e^{2\pi i b})(1-e^{2\pi i (c-b)})}\int_{\gamma}t^{b-1}(1-t)^{c-b-1}(1-tu)^{-a}dt,$$ where $\Gamma$ is the Euler Gamma-function and $\gamma$ is the Pochhammer contour. Thus, if put $(z_2,z_3,z_4)=(0,1/2,1)$, and $(a,b,c|u)=(1/3,-1/3,1/3|2)$, we see that $a_{1,2}$ is given up to a constant by ${}_2F_{1}(1/3,-1/3,1/3|2)$ which can easily be checked to be non-zero using Wolfram Mathematica. This gives a contradiction.

\end{proof}

\begin{remark}Proposition \ref{nonsplit} therefore implies that the variation of Hodge structures induced on $\mathbb{A}(\vec{z})$ by Theorem \ref{secondmain} is not always pure. Now if we consider the case $\ell=2$, i.e. $\kappa=4$, the situation is different. We know that the dimension of the conformal block $\mathbb{V}_{\mathfrak{sl}(2),\vec{\lambda},2}(\mathbb{P}^1,z_1,\dots,z_n)$ is two and is equal to the space $\mathbb{A}(\vec{\lambda})$. Since the KZ connection on conformal blocks is unitary, it follows that the monodromy representation of the KZ connection on $\mathbb{A}(\vec{\lambda})$ is semisimple for $\kappa=4$.
\end{remark}

\subsection{The weight filtration}
We examine the mixed Hodge structure in the above example (with $\ell=1$). Here $W=\Bbb{C}^2$ with coordinates $t_1$ and $t_2$. The projective compactification $P$ can be taken to be $\Bbb{P}^1\times \Bbb{P}^1$ with $(\infty,\infty)$ and the points $(z_i,z_i)$ blown up. By an easy calculation,
\begin{itemize}
\item
The exceptional divisor $E_{\infty}$ at the blow up at $(\infty,\infty)$  has weight $a_{\alpha}$ equal to $-2$.
\item
The exceptional divisors for the blow up at $(z_i,z_i)$  have weight zero,
\item The rest of the boundary divisors $E_{\alpha}$ have non-integral weights.
\end{itemize}
Let $\ml=\ml(a)$ be the local system over $U$. Consider the composite map
\begin{equation}\label{map1}
\mathbb{A}(\vec{\lambda})^*\to H^2(P,\underline{\ml}(a))\to H^2(U,\ml(a)).
\end{equation}
Let $V\subset P$ be the complement of all boundary divisors but $E_{\infty}$. Let
$\hat{U}$ be as before (Section \ref{honda}). It is easy to see that $\hat{U}\to U$ extends to an etale covering $\hat{V}\to V$. Let $\hat{Z}=\hat{V}-\hat{U}$.

Now $H^2(U,\mathcal{L}(a))$ is an isotypical component of $H^2(\hat{U},\Bbb{C})$, which maps to  $H^1(\hat{Z},\Bbb{C})$ by the Gysin map (a residue).  Therefore we obtain a map
\begin{equation}\label{map2}
\mathbb{A}(\vec{\lambda})^*\to H^1(\hat{Z},\Bbb{C}).
\end{equation}
\begin{enumerate}
\item The line of dual conformal blocks (alt level $1$) in $\mathbb{A}(\vec{\lambda})^*$ maps to zero under \eqref{map2}, but to a non zero element under \eqref{map1}(this is clear because conformal blocks give forms that extend across $\hat{Z}$ and hence have zero residues).
\item The image of any element in $\mathbb{A}(\vec{\lambda})^*$ in $H^1(\hat{Z},\Bbb{C})$ comes from the $(1,0)$ part of the smooth projective compactification of  $\hat{Z}$. To do this we can take the element of $\mathbb{A}(\vec{\lambda})^*$ to be the element in $V_{\omega_1}^{\tensor 4}\to\Bbb{C}$ which sends $a\tensor b\tensor c\tensor d$ to $(a,b)(c,d)$ where $(a,b)=a\wedge b\in \wedge^2 V_{\omega_1}=\Bbb{C}$. The assertion then reduces to: $(x(x-1))^{-\frac{2}{3}}dx$ is square integrable on $\Bbb{A}^1-\{0,1\}$, which is immediate.
\end{enumerate}
Therefore the weights filtration $\mathbb{A}(\vec{\lambda})^*$ is of the form $W_2\subset W_3$.

\subsubsection{}In the same example with $\ell=0$, one can see that the corresponding MHS is pure of weight $3$, and is (projectively) a Tate twist of the first cohomology of the elliptic curve corresponding to the point $(z_1,z_2,z_3,z_4)\in \overline{M}_{0,4}=\Bbb{P}^1$ (the Legendre family).

\subsubsection{} Let $\lambda_1,\dots,\lambda_n$ be weights for $\mathfrak{sl}(2)$
with $\sum\lambda_i=2m$, assume $0<\lambda_i<m$. Let $\mu_i=\frac{\lambda_i}{m}$ so that $\sum \mu_i=2$.
\begin{question}
Is the KZ local system, at level $\ell=m-2$, with fibers
$\mathbb{A}(\vec{\lambda})^*/\Bbb{V}^*_{\frg,\vec{\lambda},\ell}(\Bbb{P}^1,z_1,\dots,z_n)$
over the configuration space of $n$-points on $\Bbb{A}^1$ the same, up to a Tate twist,  as the one considered by Deligne and Mostow \cite{DM} (i.e., the Lauricella system) corresponding to the weights $\mu_i$?
\end{question}
\subsection{Acknowledgements}
We thank N. Fakhruddin, G. Pearlstein, E. Looijenga, A. Nair, M. Nori and A. Varchenko for useful discussions and communication. P. Brosnan and S. Mukhopadhyay were partially supported by NSF grant DMS 1361159 (PI: P. Brosnan). S. Mukhopadhyay was also supported by a Simons Travel Grant.

\noindent
P.~Belkale: Department of Mathematics, University of North Carolina, Chapel Hill, NC 27599, USA\\
{{email: belkale@email.unc.edu}}

\vspace{0.08 cm}

\noindent
P. ~Brosnan: Department of Mathematics, University of Maryland, College Park, MD 20742, USA\\
{{email: pbrosnan@umd.edu}}

\vspace{0.08 cm}

\noindent
S. ~Mukhopadhyay: Department of Mathematics, University of Maryland, College Park, MD 20742, USA\\
{{email:  swarnava@umd.edu}}

\end{document}